\documentclass[a4paper,12pt,oneside,reqno]{amsart}
\usepackage{hyperref}
\usepackage[headinclude,DIV13]{typearea}
\areaset{15.1cm}{25.0cm}
\usepackage{amsfonts,amssymb,amsmath,amsthm}
\usepackage[utf8]{inputenc}
\usepackage{graphicx}

\newtheorem{theorem}{Theorem}[section]
\newtheorem{lemma}[theorem]{Lemma}
\newtheorem{proposition}[theorem]{Proposition}

\theoremstyle{definition}

\theoremstyle{remark}
\newtheorem*{remark}{Remark}

\numberwithin{equation}{section}

\def\dist{\mathop{\rm dist}\nolimits}
\def\Var{\mathop{\rm Var}\nolimits}
\def\Cov{\mathop{\rm Cov}\nolimits}

\def\d{\mathrm{d}}

\def\bbone{\boldsymbol 1 }

\def\<{\langle}
\def\>{\rangle}
\def\a{\alpha}
\def\b{\beta}
\def\e{\epsilon}

\def\g{\gamma}

\def\s{\sigma}
\def\t{\tau}

\def\O{\Omega}

\begin{document}
\title{Universality of the REM for dynamics of mean-field spin glasses} 

\author[G. Ben Arous]{G\'erard Ben Arous}
\address{G. Ben Arous\\
  Courant Institute of the Mathematical Sciences\\
  New York University\\
  251 Mercer Street\\ 
  New York, NY 10012, USA}
\email{benarous@cims.nyu.edu}

\author[A. Bovier]{Anton Bovier}
\address{A. Bovier\\
  Weierstrass Institute for Applied Analysis and Stochastics\\
  Mohrenstrasse 39\\ 
  10117 Berlin, Germany
  \\and\\
  Mathematics Institute\\
  Berlin University of Technology\\
  Strasse des 17. Juni 136\\ 
  10269 Berlin, Germany}
\email{bovier@wias-berlin.de}

\author[J. \v Cern\'y]{Ji\v r\'\i~\v Cern\'y}
\address{J. \v Cern\'y\\
  \'Ecole Polytechnique F\'ed\'erale de Lausanne\\
  1015 Lausanne\\
  Switzerland }
\email{jiri.cerny@epfl.ch}

\subjclass[2000]{82C44,60K35,60G70}
\keywords{aging, universality, spin glasses, SK model, random walk}

\date{\today}
\begin{abstract}
  We consider a  version of a Glauber dynamics for a $p$-spin 
  Sherrington--Kirkpatrick model of a spin glass that can be seen as a 
  time change of simple random walk on the $N$-dimensional hypercube. We 
  show that, for any $p \geq  3$ and any inverse temperature $\beta>0$, 
  there exist constants $\g_0>0$, such that for all exponential time 
  scales,  $\exp(\gamma N)$, with $\g\leq \g_0$, the properly rescaled 
  \emph{clock process} (time-change process), converges to an $\a$-stable 
  subordinator where ${\a=\g/\b^2<1}$. Moreover,  the dynamics exhibits 
  aging at these time scales with time-time correlation function 
  converging to the arcsine law of this \hbox{$\alpha$-stable} 
  subordinator. In other words, up to rescaling, on these time scales 
  (that are shorter than the equilibration time of the system), the 
  dynamics of $p$-spin models ages in the same way as the REM, and by 
  extension Bouchaud's REM-like trap model, confirming the latter as a 
  universal aging mechanism for a wide range of systems. The SK model (the 
    case~$p=2$) seems to belong to a different universality class.
\end{abstract} 

\maketitle

\section{Introduction and results}
\label{s:intro}

Aging has become one of the main paradigms to describe the long-time 
behavior of complex and/or disordered systems. Systems that have strongly 
motivated this research are \emph{spin glasses}, where aging was first 
observed experimentally in the anomalous relaxation patterns of the 
magnetization \cite{LSNB83,Cha84}. The theoretical modeling of aging phenomena 
took a major leap with the introduction of so-called \emph{trap models} by 
Bouchaud and Dean in the early 1990'ies \cite{Bou92,BD95} (see \cite{BCKM98} for 
  a review). These models reproduce the characteristic power law behavior 
seen experimentally while being sufficiently simple to allow for detailed 
analytical treatment. While trap models are heuristically motivated to 
capture the behavior of the dynamics of spin glass models, there is no 
clear theoretical, let alone mathematical derivation of these from an 
underlying spin-glass dynamics. The first attempt to establish such a 
connection was made in \cite{BBG02,BBG03,BBG03b} where it was shown that 
starting from a particular Glauber dynamics of the Random Energy Model
(REM), at 
low temperatures and at the time scale slightly shorter than the equilibration 
time of the dynamics, the aging of the time-time correlation function of 
the dynamics converged to that given by Bouchaud's REM-like trap model.

On the other hand, in a series of papers \cite{BC05,BCM06,BC06b,BC07} a 
systematic investigation of a variety of trap models was initiated. In this 
process,  it emerged that there appears to be an almost 
universal aging mechanism based on $\a$-stable subordinators that governs 
aging in most of the trap models. It was also shown that the same feature 
holds for the dynamics of the REM at shorter time scales than those considered 
in \cite{BBG03,BBG03b}, and that this also happens at high temperature provided
appropriate time scales are considered \cite{BC06b}. For a general review
on trap models see \cite{BC06}.

In all models considered so far, however, the random variables describing the 
quenched disorder were considered to be independent, be it in the REM or in the
trap models. Aging in correlated spin glass models was investigated rigorously
only in some cases of spherical SK models and at very short time scales 
\cite{BDG01}. In the present paper we show for the first time that 
the same type of aging mechanism is relevant also in 
correlated spin glasses, 
at least on time scales that are short compared to equilibration time (but 
exponentially large in the volume of the system).

Let us first describe the class of models we are considering. Our state spaces will be the $N$-dimensional hypercube,
 $\mathcal S_N\equiv \{-1,1\}^N$.  
$R_N: \mathcal S_N\times \mathcal S_N\to [-1,1]$ denotes as  usual the 
normalized
overlap, $R_N(\sigma ,\tau )\equiv N^{-1}\sum_{i=1}^N \sigma_i \tau_i$. The 
Hamiltonian  of the $p$-spin SK-model is defined as $\sqrt N H_N$, where 
$H_N: \mathcal S_N \to \mathbb R$ is  the 
 centered normal process indexed
by $\mathcal S_N$ with covariance 
\begin{equation}
  \mathbb E[H_N(\sigma )H_N(\tau )]=R_N(\sigma ,\tau )^p,
\end{equation}
and $p\in \mathbb N$, $p>2$.  
We will denote by $\mathcal H$ the $\s$-algebra generated by the random variables 
$H_N(\s), \s\in \mathcal S_N,N\in \mathbb N$. 
The corresponding Gibbs measure is then given by
\begin{equation}\label{gibbs} 
\mu_{\b,N}(\s)\equiv Z^{-1}_{\b,N}e^{\b\sqrt NH_N(\s)},
\end{equation}
where $Z_{\b,N}$ denotes the normalizing partition function.

We define the classical trap-model dynamics as a nearest neighbor
continuous time Markov chain $\sigma_N(\cdot)$ on $\mathcal S_N$ with
transition rates
\begin{equation}
  \label{e:rates}
  w_N(\s,\t) =
  \begin{cases}
    N^{-1}e^{-\b \sqrt NH_N(\s)},
    &\text{if } \dist(\sigma ,\tau )=1,\\
    0,&\text{otherwise;}
  \end{cases}
\end{equation}
here $\dist(\cdot,\cdot)$ is the graph distance on the hypercube,
\begin{equation}
  \dist(\sigma ,\tau )=\frac 12 \sum_{i=1}^N |\sigma_i-\tau_i|.
\end{equation}
A simple way to construct this dynamics is as a time change of a simple
random walk on $\mathcal S_N$: We
denote by $Y_N(k)\in \mathcal S_N$, $k\in \mathbb N$, the simple 
unbiased random walk (SRW) on $\mathcal S_N$ started at some fixed point of
$\mathcal S_N$, say at $\{1,\dots,1\}$.
For $\beta>0$ we define the 
\textit{clock-process} by
\begin{equation}
  S_N(k)=
  \sum_{i=0}^{k-1} e_i \exp\big\{\beta \sqrt N H_N\big(Y_N(i)\big)\big\} ,
\end{equation}
where $\{e_i,i\in \mathbb N\}$ is a sequence of mean-one i.i.d.~exponential
random variables. We  denote by $\mathcal Y$  the $\s$-algebra
generated by the SRW random  
variables $Y_N(k)$, $k\in\mathbb N$, $N\in\mathbb N$.   The $\sigma$-algebra generated by the random variables 
$e_i, i\in \mathbb N$ will be denoted by $\mathcal E$.
Then the process $\sigma_N(\cdot)$ can be written as
\begin{equation}\label{eq:glauber.1}
  \sigma_N(t)\equiv Y_N(S_N^{-1}(t)).
\end{equation}
Obviously, $\sigma_N$ 
is reversible with respect 
to the measure $\mu_{\b,N}$.
We will consider all random processes to be defined on an abstract probability 
space $(\O,\mathcal F,\mathbb P)$. Note that the three $\s$-algebras $\mathcal H$,
$\mathcal Y$, and $\mathcal E$  
are all independent under $\mathbb P$.

We will systematically use the definition of the dynamics given by 
\eqref{e:rates} or \eqref{eq:glauber.1}. This is the same as was used in the analysis of
the REM and in most work on trap models. It differs substantially
from more popular dynamics such as the Metropolis or the heat-bath
algorithm. The main difference is that in these dynamics the
trajectories are not independent of the environment and are biased
against going up in energy. This may have a substantial effect on the
dynamics, and we do not know whether our results will apply also (with
some modifications) in these cases. The fact is that we currently do
not have the tools to analyze these dynamics even in the case of the REM!

\medskip

Let $V_\alpha(t) $ be the  $\alpha $-stable subordinator with the Laplace 
transform given by 
\begin{equation}
  \mathbb E[e^{-\lambda V_\alpha (t)}]=\exp(-t\lambda^\alpha ).
\end{equation}
The main technical result on the dynamics will be the following 
theorem that provides the asymptotic behavior of the clock process.
\begin{theorem}
  \label{t:main}
  There exists a function $\zeta (p)$ such that 
  for all $p\ge 3$ and $\gamma $ satisfying 
  \begin{equation}
   0 < \gamma < \min\left(\beta^2,\zeta (p) \beta\right),
  \end{equation}
  under the conditional
  distribution $\mathbb P[\cdot|\mathcal Y]$  the law of
  the stochastic process
  \begin{equation}
    \label{e:barSdef}
    \bar S_N(t)=
    e^{-\gamma N}
    S_N\big(\big\lfloor t N^{1/2}e^{N\gamma^2/2\beta^2}\big\rfloor\big),
    \qquad t\ge 0,
  \end{equation}
  defined on the the space of càdlàg functions equipped with the  
  Skorokhod $M_1$-topology, converges, $\mathcal Y$-a.s.,  to the law of 
  $\gamma /\beta^2$-stable subordinator 
  $V_{\gamma/\beta^2}(Kt), t\ge 0$, where $K$ is a positive constant
  depending on $\gamma $, $\beta $ and $p$.  
  
  Moreover, the function $\zeta (p)$ is increasing and it satisfies 
  \begin{equation}
    \label{e:zeta}
    \zeta (3)\simeq 1.0291
    \qquad \text{and} \qquad
    \lim_{p\to \infty}\zeta (p)=\sqrt{2\log 2}.
  \end{equation}
\end{theorem}

We will explain in Section~\ref{s:clock} what the $M_1$-topology is. 
Roughly, it is a weak topology that does not convey much information at 
the jumps of the limiting process: it can be the case that the 
approximating processes jumps several times at rather short distances to 
produce one bigger jump of the limit process. This will actually be the 
case in our models for $p<\infty$, while it is not the case in the REM. 
Therefore we cannot replace the $M_1$ topology with the stronger $J_1$-topology 
in Theorem \ref{t:main}. 

To control the behavior of spin-spin correlation functions that are 
commonly used to characterize aging, we need to know more on how these 
jumps occur at finite $N$. What we will show, is that if we the slightly 
coarse-grain the process $\bar S_N$ over blocks of size $o(N)$, the 
rescaled process does converge in the $J_1$-topology. What this says, is 
that the jumps of the limiting  process are compounded by smaller jumps 
that are made over $\leq o(N)$ steps of the SRW. In other words, the jumps 
of the limiting process come  from waiting times accumulated in one 
slightly extended trap, and during this entire time only a negligible 
fraction of the spins are flipped. That will imply the following aging 
result. 
\begin{theorem}
  \label{t:aging}
  Let $A_N^\varepsilon (t,s)$ be the event defined  by
  \begin{equation}
    A_N^\varepsilon (t,s)=
    \{
      R_N\big(
        \sigma_N\big(te^{\gamma N}\big),
        \sigma_N\big((t+s)e^{\gamma N}\big)\big)
      \ge 1-\varepsilon\big\}.
  \end{equation}
  Then, under the hypothesis of Theorem~\ref{t:main}, for all 
  $\varepsilon \in (0,1)$, $t>0$ and $s>0$,
  \begin{equation}
    \lim_{N\to \infty }\mathbb P[
      A_N^\varepsilon (t,s)]=
    \frac{\sin\alpha \pi }{\pi }
    \int_0^{t/(t+s )} u^{\alpha -1}(1-u)^{-\alpha }\,\d u.
  \end{equation}
\end{theorem}

\begin{remark} We will in fact prove the stronger statement that
  aging in the above sense occurs along almost every random walk
  trajectory, that is
  \begin{equation}
    \label{e:aging}
    \lim_{N\to \infty }\mathbb P[
      A_N^\varepsilon (t,s)
      |\mathcal Y]=
    \frac{\sin\alpha \pi }{\pi }
    \int_0^{t/(t+s )} u^{\alpha -1}(1-u)^{-\alpha }\,\d u,
    \qquad \text{$\mathcal Y$-a.s.}
  \end{equation}
\end{remark}

Let us discuss the meaning of these results. $e^{\gamma N}$ is the time-scale
at which we want to observe the process. According to Theorem~\ref{t:main}, at
this time the random walk will make of the order of 
$ N^{1/2} e^{N\g^2/2\b^2}\ll e^{\g N}$
steps. Since this number is also much smaller than $2^N$ (as follows from
  \eqref{e:zeta}), the random
walk will essentially visit that number of sites.

If the random process $H_N$ was i.i.d.,~then the maximum of $H_N$ along
the trajectory would be 
$\big(2\ln (N^{1/2} e^{N\g^2/2\b^2})\big)^{1/2} \sim N^{1/2}\gamma /\b$, 
and the time spent in that
site would be of order $e^{\g N}$. 
Since Theorem~\ref{t:main} holds also in the i.i.d.~case, that is in the
REM (see \cite{BC06b}),
the time
spent in the maximum is comparable to the total time and
 the convergence to the $\alpha $-stable subordinator
implies that the total
accumulated time is composed of pieces of order $e^{\gamma N}$ that are
collected along the trajectory. In fact, each jump of the
subordinator corresponds to one visit to a site that has waiting
times of that order. In a common metaphor, the sites are referred to as
traps and the mean waiting times as their depths.  
  
The theorem in the general case states that in the $p$-spin model, the
same is essentially true. The difference will be that the traps
here will not consist of a single site, but consist of a deep valley
(along the trajectory) whose bottom that has approximately 
the same energy as in the i.i.d.~case and whose
shape and width we will be able to describe quite
precisely. Remarkably, the number of sites contributing
significantly to the residence time in the valley is essentially
finite, and different valleys are statistically independent. 

The fact that traps are finite may appear quite surprising to those
familiar with the statics of $p$-spin models. From the results there
(see \cite{Tal03,Bov06}), it is known that the Gibbs measure
concentrates on ``lumps'' whose diameter is of order $N\e_p$, with
$\e_p>0$. 
The mystery is however solved easily: the process $H_N(\s)$
does indeed decreases essentially linearly with speed $ N^{-1/2}$  from a local
maximum.  Thus, the residence times in such sites decrease
geometrically, so that the contributions of a neighborhood of size $K$
of a local maximum amounts to a fraction of $(1-c^{-K})$ of the total
time spend in that valley ; for the support of the Gibbs measure, one
needs however to 
take into account the entropy, that is
that the volumes of the balls of radius $r$ increases like 
$N^r$.
For the dynamics, at least at our
time-scales, this is, however, irrelevant, since the SRW leaves a
local minimum essentially ballistically.

The proof of Theorem \ref{t:main} relies on the combination of detailed 
information on the properties of simple random walk on the hypercube, 
which is provided in Section~\ref{s:rwprop} (but see also \cite{Mat89,BG06,CG06}), and
comparison of the process $H_N$ on the trajectory of the SRW to a simpler
Gaussian process using interpolation techniques à la Slepian,
familiar from extreme value theory of Gaussian processes.

Let us explain this in more detail. On the time scales we are considering, the 
SRW makes 
$t N^{1/2}\exp(N\g^2/2\b^2)\ll t N^{1/2}\exp(N\zeta(p)^2/2) \ll 2^N$ 
steps.
In this regime the SRW is extremely ``transient'', in the sense that (i) 
starting from a given point $x$, for a times 
$t\leq \nu \sim N^\omega, \omega< 1$, 
the 
distance from $x$ grows essentially linearly with speed one, that is there 
are no backtrackings with high probability; 
(ii) 
the SRW will \emph{never} return to a neighborhood of size $\nu$ of the 
starting point $x$, with high probability. The upshot is that we can think of 
the trajectory of the SRW essentially as of a straight line.

Next we consider the Gaussian process restricted to the SRW trajectory.
We expect that the main contributions to the sums $S_N(k)$ come from places 
where $Y_N$ is maximal (on the trajectory). We expect that the distribution of
these extremes do not feel the correlation between points  farther 
than $\nu$ apart. On the other hand, for points closer than $\nu$, the 
correlation function $R_N(Y_N(i),Y_N(j))^p$
 can be well approximated by a linear function 
$1-2p|i-j|/N$ (using that $R_N(Y_N(i),Y_N(j))\sim 1-2|i-j|/N$). This is convenient
 since this process has an explicit representation in terms of i.i.d.~random 
variables that allow for explicit computations (in fact, this is one of the 
famous Slepian processes for which the extremal distribution can be computed 
explicitly \cite{Sle61,She71}). Thus the idea is to cut the SRW trajectory into
blocks of length $\nu$ and to replace the original process $H_N(Y_N(i))$ 
by a new one $U_i$, where $U_i$ and $U_j$ are independent, if $i,j$ are not in 
the same block, and $\mathbb E[U_iU_j]=1-2p|i-j|/N$ if they are. For the new process, 
Theorem \ref{t:main} is relatively straightforward. The main step is the
computation of Laplace transforms in Section~\ref{s:oneblock}. Comparing the real process 
with the auxiliary one is the bulk of the work and is done in
Section~\ref{s:compar}. The 
properties of SRW needed are established in Section~\ref{s:rwprop}.
In Section~\ref{s:clock} we present the proofs of the main theorems.

Our results here show some universality of the REM for dynamics of $p$-spin 
models with $p\ge 3$. This dynamic universality  is close to the static 
universality of the REM, which shows that various features of the 
landscape of energies (that is of the Hamiltonian $H_N$) are insensitive 
to correlations. This static universality in a microcanonical context has 
been introduced by \cite{BM04} (see \cite{BK06,BK06b} for rigorous results 
  on spin-glasses). The static results closest to our dynamics question are 
given in \cite{BGK06,BK07} where it is shown that the statistics of 
extreme values for the restriction of $H_N$ to a random sets 
$X_N\subset \mathcal S_N$ are universal, for $p\ge 3$ and $|X_N|=e^{cN}$, 
for $c$ small enough.

\thanks This work was initiated during a concentration period on 
metastability and aging at the Max-Planck Institute for Mathematics in the 
Sciences in Leipzig. GBA and AB thank the MIP-MIS and Stefan Müller for 
kind hospitality during this event. AB's research  is  supported in part 
by DFG in the Dutch-German Bilateral Research Group ``Mathematics of 
Random Spatial Models from Physics and Biology''.

\section{Behavior the one-block sums} 
\label{s:oneblock}

In this section we analyze the distribution of the block-sums 
$\sum_{i=1}^\nu e_i  e^{\beta \sqrt N U_i}$, 
where $e_i$ are mean-one i.i.d.~exponential  random variables, 
and $\{U_i,i=1,\dots,\nu \}$ is a centered Gaussian 
process with the covariance $\mathbb E U_i U_j= 1-2p|i-j|/N$; $\nu = \nu_N$ 
is a function of $N$ of the form 
\begin{equation}
  \label{e:nudef}
  \nu = \lfloor N^{\omega }\rfloor, \qquad \text{with }
  \omega \in (1/2,1).
\end{equation}
As explained in the introduction,
this process  will serve as a local approximation of the corresponding 
block sums along a SRW trajectory. 
We characterize the distribution of the block-sums in terms of its 
Laplace transform
\begin{equation}
  \label{e:mathcalf}
  \mathcal F_N(u ) =
  \mathbb E\Big[\exp \Big\{
      -ue^{-\gamma N} \sum_{i=1}^\nu e_i e^{\beta \sqrt N U_i}
      \Big\}\Big].
\end{equation}

\begin{proposition}
  \label{p:oneblock}
  For all $\gamma $ such that 
  $\gamma/ \beta^2\in  (0,1)$ there exists a constant,
  $K=K(\gamma ,\beta ,\omega ,p )$,  such that, uniformly for $u$ in compact
  subsets of $[0,\infty)$,
  \begin{equation}
    \lim_{N\to\infty}
    N^{1/2}\nu^{-1} 
    e^{N \gamma^2/2\beta^2}\,
    [ 1-\mathcal F_N(u)]=
    K u^{\gamma /\beta^2}.
  \end{equation}
\end{proposition}
\begin{proof}
  We first compute the conditional expectation in \eqref{e:mathcalf} given 
  the $\s$-algebra, $\mathcal U$, generated by the Gaussian process $U$, 
  \begin{equation}
    \label{e:ab}
    \begin{split}
      \mathbb E\Big[\exp \Big\{
          -ue^{-\gamma N} \sum_{i=1}^\nu e_i e^{\beta \sqrt N U_i}
          \Big\}
        \Big|\mathcal U\Big]
      &=\prod_{i=1}^\nu \frac 1{1+u e^{-\gamma N}  e^{\beta \sqrt N U_i}}
      \\&=
      \exp\bigg\{-\sum_{i=1}^\nu g\Big(u e^{-\gamma N} 
          e^{\beta \sqrt N U_i}\Big)\bigg\},
    \end{split}
  \end{equation}
  where 
  \begin{equation}
    \label{e:g}
    g(x)\equiv \ln(1+x).
  \end{equation}
  Note that importantly, $g(x)$ is monotone increasing and non-negative
  for $x\in\mathbb R_+$.
  We use the well-known fact (see e.g.~\cite{Sle61}) 
that the random variables 
  $U_i$ can be expressed using a sequence of i.i.d.~standard normal 
  variables, $Z_i$, as follows. Set  
  $Z_1=(U_1+U_\nu )/(4-4p(\nu -1)/N)^{1/2}$
  and $Z_k=(U_k-U_{k-1})/(4p/N)^{1/2}$, $k=2,\dots, \nu $. Then  
  $Z_i$ are i.i.d.~standard normal and 
  \begin{equation}
    U_i=\Gamma_1 Z_1 + \dots + \Gamma_i Z_i 
    - \Gamma_{i+1}Z_{i+1} - \Gamma_\nu  Z_\nu ,
  \end{equation}
  where 
  \begin{equation}
    \Gamma_1=\sqrt{1-\frac p{N}(\nu -1)}\qquad
    \text{and}\qquad
    \Gamma_2=\dots=\Gamma_\nu =\sqrt{\frac p{N}}.
  \end{equation}
  Observe that $\sum_{i=1}^\nu \Gamma_i^2 = 1$.
  Let us define $G_i(z) = G_i(z_1,\dots,z_\nu )$ as
  \begin{equation}
    G_i(z)=\Gamma_1 z_1 + \dots + \Gamma_i z_i
    -\Gamma_{i+1}z_{i+1}-\dots-\Gamma_\nu  z_\nu .
  \end{equation}
  Using this notation we get 
  \begin{equation}
    \label{e:aa}
    1-\mathcal F_N(u )
    =
    \int_{\mathbb R^\nu }
    \frac {\d z}{(2\pi )^{\nu /2}}
    e^{-\frac 12\sum_{i=1}^\nu z_i^2 }
    \Big\{ 1-
      \exp\Big[
        -\sum_{i=1}^\nu g\left(u e^{-\gamma N} e^{\beta \sqrt N
            G_i(z)}\right)\Big]\Big\}.
  \end{equation}
  We divide the domain of integration into several parts according to which
  of the  $G_i(z)$ is 
  maximal. Define $D_k=\{z:G_k(z)\ge G_i(z)\forall i\neq k\}$. On 
  $D_k$ we use the substitution 
  \begin{equation}
    \label{e:bsubst}
    \begin{split}
      z_i&=b_i+\Gamma_i(\gamma N - \log u)/(\beta \sqrt N),
      \qquad\text{if $i\le k$,}\\
      z_i&=b_i-\Gamma_i(\gamma N - \log u)/(\beta \sqrt N),
      \qquad\text{if $i> k$.}
    \end{split}
  \end{equation}
  It will be useful to define $\sum_{j=i+1}^k a_j$ as 
  $\sum_{j=1}^k a_j - \sum_{j=1}^i a_j$, which is meaningful also for $k<i+1$.
  Using this definition 
  \begin{equation}
    G_k(b)-G_i(b)=2\sum_{j=i+1}^k \Gamma_\nu b_j.
  \end{equation}
  Set $\theta =-\log(u)/(\gamma N)$ and define 
  \begin{equation}
    D'_k=\Big\{b:\sum_{j=i+1}^k b_j+\frac {\gamma\sqrt p } \beta 
      |k-i|(1+\theta )\ge 0\,\forall i\neq k\Big\}.
  \end{equation}
  After a straightforward computation we find that \eqref{e:aa} equals
  \begin{equation}
    \begin{split}
      \label{e:bb}
      e^{-N{\gamma^2}/{2\beta^2}}&
      u^{\gamma /\beta^2}
      \sum_{k=1}^\nu 
      \int_{ D'_k }
      \frac {\d b}{(2\pi )^{\nu /2}}
      e^{-\frac 12\sum_{i=1}^\nu  b_i^2}
      e^{-\frac \gamma \beta \sqrt N G_k(b)(1+\theta )}
      \\&\times
      \Big\{1- 
        \exp\Big(-
          \sum_{i=1}^\nu
          g\left(e^{\beta \sqrt N G_k(b)-
            2\beta \sqrt{p}\sum_{j=i+1}^k b_j -2p\gamma|k-i|(1+\theta )
      }\right)\Big)\Big\}.
    \end{split}
  \end{equation}
  To finish the proof we have to show that $u^{\gamma /\beta^2}$ is 
  asymptotically the only dependence of \eqref{e:bb} on $u$ (or on 
  $\theta $) and that the 
  sum is of order $\nu N^{-1/2}$. 
  We change variables once more to $a_j=b_j/(1+\theta)$ in order to
  remove the dependence of the integration domains on $u$.
  Then the sum (without the prefactor) 
  in \eqref{e:bb} can be expressed as
  \begin{equation}
    \begin{split}
      \label{e:cc}
      \sum_{k=1}^\nu &
      \int_{D''_k }
      \frac {(1+\theta )^\nu \d a}{(2\pi )^{\nu /2}}
      e^{- \frac 12(1+\theta )^2\sum_{i=1}^\nu  a_i^2}
      \bigg[
      e^{-\frac \gamma \beta \sqrt N G_k(a)(1+\theta )^2}
      \\&\times
      \Big\{1- 
        \exp\Big(-\sum_{i=1}^\nu
          g\left(
            e^{(\beta \sqrt N G_k(a)- 
                2\beta \sqrt{p}\sum_{j=i+1}^k a_j 
                -2p\gamma|k-i|)
              (1+\theta )
      }\right)\Big)\Big\}\bigg],
    \end{split}
  \end{equation}
  where $D''_k=\big\{a:\sum_{j=i+1}^k a_j+\frac {\gamma\sqrt p} \beta 
    |k-i|\ge 0\,\forall i\neq k\big\}.$

  Let $\delta >0$ be such that $(1+\delta )\gamma /\beta^2 < 1$, and let 
  $N > \log(u)/(\gamma\delta)$, so  that $|\theta |\le \delta $.  
  We first examine the bracket in 
   the above expression for a fixed $k$. On $D''_k$
  \begin{equation}
    \exp\Big\{
      -\sum_{i=1}^\nu
      g\big(e^{(\beta \sqrt N G_k(a)-
          2\beta \sqrt{p}\sum_{j=i+1}^k a_j -2p\gamma|k-i|)(1+\theta )
  }\big)\Big\}\ge
    \exp\big\{-\nu g\big(e^{\beta \sqrt N G_k(a)(1+\theta )}\big)\big\}.
  \end{equation}
  Write $G_k(a)$ as  (recall \eqref{e:nudef})
  \begin{equation}
    G_k(a)=\frac {\xi -\omega \log N}{(1+\theta ) \beta \sqrt N}.
  \end{equation}
  The bracket  of \eqref{e:cc} is then smaller than
  \begin{equation}
    \label{e:ab.2}
    \begin{split}
      & e^{-\frac \gamma {\beta^2}  (\xi -\omega \log N)(1+\theta )}
      \big\{1-  
        \exp\big(-\nu g\big(e^{\xi -\omega \log N}\big)\big)
        \big\}
      \\ &\quad =
      N^{\frac{\gamma \omega  (1+\theta )}{\beta^2}}
      e^{-\frac {\gamma\xi } {\beta^2}  (1+\theta )}
      \big\{1-  
        \exp\big(- \nu g\big(e^{\xi}/\nu\big)\big)
        \big\}.
    \end{split}
  \end{equation}
  The function  
  $e^{-\frac {\gamma\xi } {\beta^2}  (1+\theta )} 
  \big\{1-  \exp\big(-  \nu g\big(e^{\xi}/\nu\big)\big) \big\}$ 
  is bounded for $\xi \in \mathbb R$, uniformly in $\nu$, 
 if $(1+\theta )\gamma /\beta^2 <1$.  
 Namely, if $\xi\ge 0 $, 
 \begin{equation}
   e^{-\frac {\gamma\xi } {\beta^2}  (1+\theta )} 
   \big\{1-  \exp\big(-  \nu g\left(e^{\xi}/\nu\right)\big) \big\}
   \leq e^{-\frac {\gamma\xi } {\beta^2}  (1+\theta )} \leq 1.
 \end{equation}
 If $\xi<0$, then, since $g(x)\le x$,
 \begin{equation}
   \big\{1-  \exp\big(-  \nu g\big(e^{\xi}/\nu\big)\big) \big\}
   \le
   \big\{1-  \exp\big(-  e^{\xi}\big) \big\},
 \end{equation}
 which behaves like $e^{\xi}$, as $\xi\to -\infty$. This  
 compensates the exponentially growing prefactor, if 
 $(1+\theta )\gamma /\beta^2 <1$. Thus, under this condition,
 the bracket of \eqref{e:cc} 
 increases at most polynomially with $N$.

  In view of this at most polynomial increase, there exist $\delta>0$ 
  small, such that  the domain of integration in \eqref{e:cc} may be 
  restricted to $a_i$'s satisfying 
  \begin{equation}
    \label{e:domain}
    \nu^{-1}\sum_{i=1}^\nu a_i^2\in (1-\delta ,1+\delta ),
    \quad |a_1|\le N^{1/4}, \quad \sum_{i=1}^\nu |a_i| \le \nu^{1+\delta}.
  \end{equation}
  The integral over the remaining $a_i$'s decays at least as 
  $e^{-N^{\delta'}}$ for some $\delta '>0$ (by a simple large deviation
    argument).  
  For all $a$ satisfying \eqref{e:domain}, 
  $|G_k(a)| \le N^{1/4}+N^{-1/2}\nu^{1+\delta'}\ll N^{1/2}$ and thus, 
  for any fixed $u$, uniformly in $a$,
  \begin{equation}
    \frac{ e^{-\frac \gamma \beta 
        \sqrt N G_k(a)(1+\theta )}}
    {e^{-\frac \gamma \beta \sqrt N G_k(a)}}
    \xrightarrow{N\to \infty} 1,
    \quad\text{and}\quad
    \frac{e^{-\frac 12 (1+\theta )^2 \sum_{i=1}^\nu  a_i^2}}
    {e^{-\frac 12 \sum_{i=1}^\nu  a_i^2}}
    \xrightarrow{N\to \infty} 1.
  \end{equation}
  Also, $ (1+\theta )^\nu \xrightarrow{N\to \infty} 1 $. Hence, up to a
  small error, we can remove all but the last occurrence of $\theta $
  in \eqref{e:cc}.

  Finally, taking $x_i=a_i$ for $i\ge 2$, $x_1=N^{1/2} G_k(a)$, and thus
  \begin{equation}
    a_1=\frac{x_1-4p(x_2+\dots+x_k-x_{k+1}-\dots-x_\nu) }{\Gamma_1\sqrt N},
  \end{equation}
  \eqref{e:cc} equals, up to a small error,
  \begin{equation}
     \begin{split}
      \label{e:dd}
      &\sum_{k=1}^\nu 
      \int_{D''_k }
      \frac {\d x\, e^{- \frac 12\sum_{i=2}^\nu  x_i^2}}
      {\Gamma_1 N^{1/2}(2\pi )^{\nu /2}}
      \exp\Big(-\frac \gamma \beta x_1-\frac {x_1^2}{2 \Gamma_1^2 N}\Big)
      \exp\Big(-\frac {a_1^2} 2+\frac {x_1^2}{2 \Gamma_1^2 N}\Big)
      \\&\times
      \Big\{1- 
        \exp\Big(-\sum_{i=1}^\nu g\Big(e^{(1+\theta )\beta x_1}
          e^{-\left(
            2\beta \sqrt{p}\sum_{j=i+1}^k x_j -2p\gamma|k-i|\right)(1+\theta )
    }\Big)\Big)\Big\}.
    \end{split}
  \end{equation}
  The last exponential term on the first line can be omitted. Indeed,
  \begin{equation}
    -\frac {a_1^2} 2+\frac {x_1^2}{2 \Gamma_1^2 N}=
    \frac 4 { \Gamma_1^2 N}
    \big[px_1(x_2+\dots-x_\nu)-2p^2(x_2+\dots-x_\nu)^2\big]
    \xrightarrow{N\to\infty}0
  \end{equation}
  uniformly for all $|x_1|\le N^{(1+\delta)/2} $ and 
  $|x_2+\dots-x_\nu|\le \nu^{(1+\delta)/2 }$, if $\delta >0$ sufficiently
  small. The integral over the remaining 
  $x$ is again at most $e^{-N^{\delta '}}$.  
  
  Now we estimate the integral
  over $x_2,\dots,x_\nu $,
  \begin{equation}
      \label{e:ee}
      \int_{\bar D''_k }
      \frac {\d x  e^{- \frac 12\sum_{i=2}^\nu  x_i^2}}
      {(2\pi )^{(\nu -1) /2}}
      \exp\Big(-\sum_{i=1}^\nu g\Big(e^{(1+\theta )\beta x_1}
        e^{-\left(
            2\beta \sqrt{p}\sum_{j=i+1}^k x_j +2p\gamma|k-i|\right)(1+\theta )}
        \Big)\Big),
  \end{equation}
  where $\bar D''_k$ is the restriction of $D''_k$ to the last $\nu -1$
  coordinates (which does not depend on the value of the first one).
  Let $V=(V_2,\dots,V_\nu) $ be a sequence of i.i.d.~standard normal random
  variables. Then,
  \eqref{e:ee} equals
  \begin{equation}\label{e:ee.2}
    \mathbb P[V\in \bar D''_k]
    \mathbb E\Big[
      \exp\Big(-\sum_{i=1}^\nu g\left(e^{(1+\theta )\beta x_1}
        e^{-\left(
            2\beta \sqrt{p}\sum_{j=i+1}^k V_j +2p\gamma|k-i|\right)
          (1+\theta )}
        \right)\Big)\Big|
       V \in \bar D_k''
      \Big].
  \end{equation}
  The probability $\mathbb P[V\in \bar D''_k]$ is bounded from below  by 
  the
  probability that the two-sided random walk, 
  $R_i=\sum_{j=0}^i V_j$, 
  $i\in \mathbb Z$,  with standard
  normal increments is larger than 
  $-{\gamma \sqrt p } |i|/\beta $ for all $i$. This probability
  is positive and does not depend on $N$, which implies that, for all $k$, 
  \begin{equation}
    \label{e:iii}
    1> \mathbb P[V\in \bar D''_k]\ge c>0.
  \end{equation}
  
  The expectation in \eqref{e:ee.2} is bounded by one, since the
  functions $g$ is positive on the domain of integration.
  Moreover, as $x_1\to-\infty$, the argument of $g$ in \eqref{e:ee.2}
  tends to zero (since the first exponential does, and the second is
    bounded by one on $D''_k$). Hence 
  \begin{equation}
    g\left(e^{(1+\theta )\beta x_1}
      e^{-(
          2\beta \sqrt{p}\sum_{j=i+1}^k V_j +2p\gamma|k-i|)(1+\theta )}
    \right)
    \sim 
    e^{(1+\theta )\beta x_1}
    e^{-(
        2\beta \sqrt{2}\sum_{j=i+1}^k V_j
        +2p\gamma|k-i|)(1+\theta )}.
  \end{equation}
  Therefore, as $x_i\to -\infty$,    
  \begin{equation}
    \begin{split}
      \label{e:ff}
      & \mathbb E\Big[
        \exp\Big(-\sum_{i=1}^\nu g\left(e^{(1+\theta )\beta x_1}
            e^{-\left(
                2\beta \sqrt{p}\sum_{j=i+1}^k V_j +2p\gamma|k-i|\right)(1+\theta )}
        \right)\Big)\Big|
        V \in \bar D_k''
        \Big]\\
      &\sim 1-e^{(1+\theta)\beta x_1 } 
      \mathbb E\Big[
        \sum_{i=1}^\nu 
        e^{-\left(
            2\beta \sqrt{p}\sum_{j=i+1}^k V_j +2p\gamma|k-i|\right)(1+\theta)}
        \Big| V\in D''_k\Big]
      \\&=
      1-e^{(1+\theta)\beta x_1 } 
      \sum_{i=1}^\nu 
      \mathbb E\Big[
        e^{-\left(
            2\beta \sqrt{p}R_{k-i}  +2p\gamma|k-i|\right)(1+\theta)}
        \Big| R_{k-i}\ge -\frac {\gamma \sqrt p }\beta |k-i| \Big].
    \end{split}
  \end{equation}
  Since $R_i$ is a centered normal random variable with variance $|i|$,  a 
  straightforward Gaussian calculation  implies that 
  \begin{equation}
    \mathbb E\Big[
        e^{-\left(
            2\beta \sqrt{p}R_{k-i}  +2p\gamma|k-i|\right)(1+\theta)}
        \Big| R_{k-i}\ge -\frac {\gamma \sqrt p }\beta |k-i| \Big]
      \sim \frac{C_{\beta ,\gamma ,p}}{\sqrt{|k-i|}}
      e^{-\gamma^2p|k-i|/(2\beta^2) }.
  \end{equation}
  Hence,  \eqref{e:ff} is essentially a summation of a geometrical sequence
  and therefore
  there exists constants $c_1$, $c_2$ independent of $k$, such that     
  \begin{equation}
    \label{e:iv}
    1-c_1 e^{(1+\theta)\beta x_1 }
    \le \eqref{e:ff}\le 
    1-c_2 e^{(1+\theta)\beta x_1 },
    \qquad\forall  x_1 < 0.
  \end{equation}
  Bounds \eqref{e:iii} and \eqref{e:iv} imply that \eqref{e:ee} is bounded 
  from above and from below (with different constants) by
  \begin{equation}
    \label{e:uio}
    C N^{-1/2}
    \exp\Big(-\frac \gamma \beta x_1-\frac {x_1^2}{2 \Gamma_1^2 N}\Big)
    (1 \wedge c e^{(1+\theta)\beta x_1}).
  \end{equation}
  and hence \eqref{e:dd} is bounded from above and below by 
  \begin{equation}
    \label{e:vv}
    C \nu N^{-1/2}\int_{\mathbb R} \d x_1
    \exp\Big(-\frac \gamma \beta x_1-\frac {x_1^2}{2 \Gamma_1^2 N}\Big)
    (1 \wedge c e^{(1+\theta)\beta x_1}) = C \nu N^{-1/2}.
  \end{equation}
  Moreover, \eqref{e:ee} is decreasing as function of $\min(k,\nu -k)$.
  As this minimum tends to infinity, \eqref{e:ee} 
  behaves as 
  $f(x_1)N^{-1/2}$ which is of course satisfy the bound \eqref{e:uio}.    
  Due to this convergence, the constants in the lower
  and the upper bound of  \eqref{e:vv} 
  can be made arbitrarily close. This completes the proof
  of Proposition~\ref{p:oneblock}.
\end{proof}

We close this section with a short description of the shape of the valleys 
mentioned in the introduction. First, it follows from \eqref{e:bsubst} and 
the following computations that the most important contribution to the 
Laplace transform comes from realizations for which 
$\max \{U_i:1\le i\le \nu \}\sim\gamma \sqrt N /\beta$ with an error of 
order $N^{-1/2}$. It is the ``geometrical'' sequence in \eqref{e:ff} which 
shows that only finitely many neighbors of the maximum actually contribute 
to the Laplace transform. The same can be seen, at least heuristically, 
from a simple calculation
\begin{equation}
  \mathbb E\Big[U_{k+i}\Big |U_k=\frac{\gamma }{\beta} \sqrt N \Big]=
  \frac {\gamma \sqrt N}{\beta } 
  - C_{\beta ,\gamma ,p} \frac{|i|}{\sqrt N}. 
\end{equation}
Which means that, disregarding the fluctuations, 
the energy decreases linearly with the distance from the
local maximum and thus the mean waiting times decrease exponentially.
  
\section{Comparison of the real and the block process}
\label{s:compar}

We now come to the main task, the comparison of the clock-process sums 
with those in which 
the  real Gaussian process is replaced by a simplified process. For a given 
realization, $Y_N$, of the SRW,
we set  $X_N^0(i)=H_N\big(Y_N(i)\big)$ (the dependence on $Y_N$ will be 
suppressed in the notation). Then 
 $X^0_N(i)$ is a centered Gaussian process indexed by $\mathbb N$ 
 with covariance matrix 
\begin{equation}
  \Lambda^0_{ij}=\mathbb E[X^0_N(i) X^0_N(j)]=R_N\big(Y_N(i),Y_N(j)\big)^p.
\end{equation}
Now we define the comparison process, $X^1_N(i)$, as the centered 
 Gaussian process with the
covariance matrix
\begin{equation}
  \Lambda^1_{ij}=\mathbb E[X^1_N(i) X^1_N(j)]=
  \begin{cases}
    1-2p|i-j|/N,&\text{if $\lfloor i/\nu \rfloor=\lfloor j/\nu\rfloor$,}\\
    0,&\text{otherwise.}
  \end{cases}
\end{equation}
For $h\in[0,1]$ we define the interpolating process
 $X^h_N(i)\equiv\sqrt{1-h}X^0_N(i)+\sqrt h X^1_N(i)$. 

Let $\ell\in \mathbb N$, $0=t_0<\dots<t_\ell=T$ and 
$u_1,\dots,u_\ell\in \mathbb R_+$ be fixed. 
For any Gaussian process $X$ we define a function 
$F_N(X)=F_N\big(X;\{t_i\},\{u_i\}\big)$ as
\begin{equation}
  \begin{split}
    F_N\big(X;\{t_i\},\{u_i\}\big)
    &\equiv
    \mathbb E\Big[\exp\Big(
        -\sum_{k=1}^\ell
        \frac {u_k} {e^{\gamma N}} 
        \sum_{i=t_{k-1}r(N)+1}^{t_k r(N)} e_i
        e^{\beta \sqrt N X(i)}\Big)
      \Big| \mathcal X
      \Big](X)
    \\&=
    \exp\Big(- 
      \sum_{k=1}^\ell
      \sum_{i=t_{k-1}r(N)}^{t_k r(N)-1} 
      g\Big(\frac {u_k} {e^{\gamma N}} 
        e^{\beta \sqrt N X(i)}\Big)
      \Big),
  \end{split}
\end{equation}
where $r(N)=N^{1/2}e^{N\gamma^2/2\beta^2}$.
Observe that 
$\mathbb E[F(X^0;t,u)|\mathcal Y]$ is a joint Laplace transform of the distribution of 
the properly rescaled clock process at times $t_i$. 
The following approximation is the crucial step of the
proof.

\begin{proposition}
  \label{p:comp}
  If the assumptions of Theorem~\ref{t:main} are satisfied, then for all
  sequences $\{t_i\}$ and $\{u_i\}$,
  \begin{equation}
    \lim_{N\to \infty} 
    \mathbb E\big[F_N\big(X^0_N;\{t_i\},\{u_i\}\big)\big|\mathcal Y\big]
    -\mathbb E\big[F_N\big(X^1_N;\{t_i\},\{u_i\}\big)\big]=0,\qquad
    \text{$\mathcal Y$-a.s.}
  \end{equation}
\end{proposition}

\begin{proof}
  We use the well-known interpolation formula for functionals of two
  Gaussian processes due (probably) to Slepian and Kahane 
  (see e.g. \cite{LT91} 
  \begin{equation}
    \label{e:comparison}
    \mathbb E[F_N(X^1_N)-F_N(X^0_N)|\mathcal Y]=
    \frac 12
    \int_0^1 \d h
    \sum_{\substack{i,j =1\\i\neq j}}^{t r(N)}
    (\Lambda^0_{ij}-\Lambda^1_{ij})
    \mathbb E\Big[
      \frac{\partial^2 F_N(X^h_N)}{\partial X(i)\partial X(j)}
      \Big|\mathcal Y
      \Big].
  \end{equation}
  We will show that the integral in \eqref{e:comparison} converges to $0$.

  Let $k(i)$ be defined by $t_{k(i)-1}r(N)< i \le t_{k(i)}r(N)$.
  The second derivative in \eqref{e:comparison} is equal to
  \begin{equation}
    \begin{split}
      \frac{u_{k(i)}u_{k(j)} \beta^2 N} {e^{2\gamma N}} 
      &e^{\beta \sqrt N(X^h_N(i)+X^h_N(j))}
      g'\Big(\frac {u_{k(i)}} {e^{\gamma N}} 
        e^{\beta \sqrt N X^h_N(i)}\Big)  
      g'\Big(\frac {u_{k(j)}} {e^{\gamma N}} 
        e^{\beta \sqrt N X^h_N(j)}\Big)     
      F_N(X^h_N)
      \\&\le 
      \frac{u_{k(i)} u_{k(j)} \beta^2 N}{e^{2\gamma N}} 
      e^{\beta \sqrt N(X^h_N(i)+X^h_N(j))}
      \\&\qquad\times 
      \exp\Big[
        -2g\Big(\frac {u_{k(i)}}{e^{\gamma N}}
          e^{\beta \sqrt N X^h_N(i)}\Big)
        -2g\Big(\frac {u_{k(j)}}{e^{\gamma N}}
          e^{\beta \sqrt N X^h_N(j)}\Big)
        \Big],
    \end{split}
  \end{equation}
  where we used that $g'(x)=(1+x)^{-1}=\exp(-g(x))$ (recall \eqref{e:g}), 
  and we omitted  in the summation of $F_N(X_N^h)$ all terms different 
  from $i$ and $j$. To estimate the expected value of this expression we 
  need the following technical lemma.
  \begin{lemma}
    \label{l:tech}
    Let $c\in [-1,1]$ and let
    $U_1$, $U_2$ be two standard normal variables with  the covariance
    $\mathbb E[U_1U_1]=c$ and $\lambda $ a small constant, 
    $0<\lambda <1-\gamma /\beta^2$ (which will stay fixed).  Define  
    $\Xi_N(c)=\Xi_N(c,\beta ,\gamma ,u,v)$ and 
    $\bar \Xi_N(c)=\bar \Xi_N(c,\beta ,\gamma ,u,v,\lambda )$ by
    \begin{equation}
      \Xi_N(c)=
      \frac{uv \beta^2 N}{e^{2\gamma N}}
      \mathbb E\Big[
        \exp\Big\{\beta \sqrt N(U_1+U_2)
          -2g\big(ue^{\beta \sqrt N U_1-  \gamma N}\big)
          -2g\big(ve^{\beta \sqrt N U_2-  \gamma N}\big)
          \Big\}
        \Big]
    \end{equation}
    and 
    \begin{equation}
      \label{e:xibarprop}
      \bar \Xi_N (c)=
      \begin{cases}
        \frac{C(\gamma ,\beta ,u,v,\lambda ) }
        {(1-c)^{1/2}}
        \exp\Big\{-\frac {\gamma^2N}{\beta^2(1+c)}\Big\},
        &\text{if $c>(\gamma/ \beta^{2})+\lambda-1 $,}
        \\
        C'(\gamma ,\beta ,u,v)N 
        \exp\big\{N(\beta^2(1+c)-2\gamma) \big\},
        &\text{if $c\le (\gamma/ \beta^{2})+\lambda-1 $,}
      \end{cases}
    \end{equation}
    where $C(\gamma ,\beta ,u,v,\lambda )$ and $C'(\gamma ,\beta ,u,v)$ 
    are suitably chosen constants, independent of $N$ and $c$. Then
    \begin{equation}
      \Xi_N(c)\le \bar \Xi_N(c).
    \end{equation}
  \end{lemma}
  \begin{proof}
    Define $\kappa_{\pm}=\sqrt{2(1\pm c)}$. Let $\bar U_1$, $\bar U_2$ be two
    independent standard normal variables. Then $U_1$ and $U_2$ can be
    written as
    \begin{equation}
        U_1=\frac 12(\kappa_+ \bar U_1 + \kappa_- \bar U_2),\qquad
        U_2=\frac 12(\kappa_+ \bar U_1 - \kappa_- \bar U_2).
    \end{equation}
    Hence, $U_1+U_2=\kappa_+\bar U_1 $. Using 
    $g(x)+g(y)=g(x+y+xy)\ge g(x+y)$ and 
    $ue^x+ve^{-x}\ge \min(u,v)e^{|x|}$, we get
    \begin{equation}
      \begin{split}
        g\big(u&e^{\beta \sqrt N U_1-  \gamma N}\big)+
        g\big(ve^{\beta \sqrt N U_2-  \gamma N}\big)
        \\&\ge
        g\Big(\min(u,v) \exp\Big(
            \frac{\kappa_+\beta \sqrt N \bar U_1}{2}+
            \Big|\frac{\kappa_-\beta \sqrt N \bar U_2}{2}\Big|
            -\gamma N
            \Big)
          \Big).
      \end{split}
    \end{equation}
    Denoting $\min(u,v)$ by $\bar u$, we find that 
    $\Xi_N(c)$ is bounded from above by
    \begin{equation}
      \frac{uv \beta^2 N}{e^{2\gamma N}}
      \int_{\mathbb R^2}\frac {\d y}{2\pi}
      \exp\Big\{
        -\frac{y_1^2+y_2^2}2
        +\beta \sqrt N \kappa_+y_1
        -2g\big(\bar u
        e^{\kappa_+\beta \sqrt N y_1/2 + 
          \kappa_-\beta \sqrt N |y_2|/2-\gamma N}\big)
        \Big\}.
    \end{equation}
    Substituting $z_1=y_1-\beta \sqrt N  \kappa_+$, $z_2=y_2$ we
    get
    \begin{equation}
      \begin{split}
        \label{e:tyu}
        &
        \frac{uv \beta^2 N}{e^{2\gamma N}}
        e^{\beta^2 \kappa_+^2 N /2}
        \int_{\mathbb R^2}\frac {\d z}{2\pi}
        \exp\Big(
          -\frac{z_1^2+z_2^2}2
          \Big)
        \\&\times
        \exp\Big(
          -2g\Big(\bar u\exp\Big\{
              \sqrt N\Big[\Big(\frac {\beta^2 \kappa_+^2}2-\gamma \Big)\sqrt N 
                +\frac{\beta \kappa_+}{2}z_1
                +\frac{\beta \kappa_-}{2}|z_2|
                \Big]
              \Big\}
            \Big)
          \Big).
      \end{split}
    \end{equation}
    The function $\exp(-2g(\bar u e^{\sqrt N x}))$ converges to the indicator 
    function  $\bbone_{x<0}$, as $N\to\infty$. The rôle of $x$ will be played by 
    the bracket in the expression \eqref{e:tyu}. 
    
    If this bracket remains negative for 
    $z$ close to zero, that is if $\gamma \ge -\lambda'+\beta^2 \kappa_+^2/2$ 
    (or equivalently $c\le (\gamma/ \beta^{2})+\lambda-1 $), then the integral in 
    \eqref{e:tyu} is bounded from above by $1$. This yields the claim of 
    the lemma for such $c$:  
    \begin{equation}
      \Xi_N(c)\le 
      \frac{uv \beta^2 N}{e^{2\gamma N}}
      e^{\beta^2 \kappa_+^2 N /2}=
      C'(\gamma ,\beta ,u,v)N 
      \exp\big\{N(\beta^2(1+c)-2\gamma) \big\} =\bar \Xi_N(c).
    \end{equation}

    If this is not the case, that is $\gamma < - \lambda '+\beta^2 \kappa_+^2/2$, then we
    need another substitution, 
    \begin{equation} \label{eq:substitute.10}
      \begin{split}
        z_1&=\frac 1{\sqrt N} \Big[ v_1 - \frac {\kappa_-}{\kappa_+} |v_2| 
          -N\Big(\beta \kappa_+ -\frac{2\gamma }{\beta \kappa_+}\Big)\Big],\\ 
        z_2&= \frac{v_2}{\sqrt N}.
      \end{split}
    \end{equation}
    This substitution transforms the domain where the bracket of 
    \eqref{e:tyu} is negative into the half-plain $v_1<0$: The expression 
    inside of the braces in \eqref{e:tyu} equals $\beta \kappa_+ v_1/2$. 
    Substituting \eqref{eq:substitute.10} into $(z_1^2+z^2_2)/2$ produces 
    an additional exponential prefactor 
    $\exp\big(-\frac{(\beta^2 \kappa_+^2-2\gamma )^2N}
      {2\beta^2 \kappa_+^2}\big)$. 
    Another prefactor $N^{-1}$ comes from the Jacobian.  
    The remaining terms can be bounded from above by
    \begin{equation}
      \label{e:uuuu}
      \int_{\mathbb R^2}
      \frac {\d v}{2\pi }
      \exp\Big\{\Big(\beta \kappa_+ - \frac {2\gamma }{\beta \kappa_+}\Big)
        \Big(v_1-\frac{\kappa_-}{\kappa_+}|v_2|\Big)
        -2g(\bar u e^{\beta \kappa_+/2})\Big\},
    \end{equation}
    which can be separated into a product of two integrals.
    The integration over $v_2$ gives a factor 
    \begin{equation}
      \Big(\Big(\beta \kappa_+ - \frac {2\gamma }{\beta \kappa_+}\Big)
        \frac{\kappa_-}{\kappa_+}\Big)^{-1}
      \le C(\lambda )\kappa_-^{-1}\le C(\lambda ) (1-c)^{-1/2}.
    \end{equation}
    Using properties of 
    $g$, the integrand of \eqref{e:uuuu} behaves as 
    $\exp\{-2v_1\gamma /\beta \kappa_+\}$ as $v_1\to\infty$, and as
    $\exp\{(\beta \kappa_+ -(2\gamma /\beta \kappa_+))v_1\}$ as 
    $v_1\to-\infty$. Therefore, the integral over $v_1$ is bounded
    uniformly by some $\lambda $-dependent constant  
    for all values of $c\ge -1+(\gamma/ \beta^{2})+\lambda $.
    Putting everything together
    \begin{equation}
      \begin{split}
        \Xi_N(c)&\le C (1-c)^{-1/2}
        \frac{uv \beta^2 N}{e^{2\gamma N}}
        e^{\beta^2 \kappa_+^2 N /2}\frac 1 N
        \exp\Big(-
          \frac{(\beta^2 \kappa_+^2-2\gamma )^2N}
          {2\beta^2 \kappa_+^2}\Big)
        \\& =
        C(\gamma ,\beta ,u,v,\lambda ) (1-c)^{-1/2}
        \exp\Big\{-\frac {\gamma^2N}{\beta^2(1+c)}\Big\}=\bar\Xi_N(c).
      \end{split}
    \end{equation}
    This finishes the proof of Lemma~\ref{l:tech}.
  \end{proof}

  Let $\|d\|=\min(d,N-d)$ and $D_{ij}=\dist(Y_N(i),Y_N(j))$.   
  Define, with a slight abuse of notation, 
  $\Lambda^0_d=(1-2dN^{-1})^p$. That is $\Lambda^0_d$ is the covariance
  of $X_N^0(i)$ and $X_N^0(j)$ if $D_{ij}=d$.
  The next proposition, which will be proved in Section~\ref{s:rwprop},
  will be used to control the correlations of the process 
  $X^0_N$. 

  \begin{proposition}
    \label{p:rwp}
    Let $\gamma $ and $\beta $ satisfy the hypothesis of 
    Theorem~\ref{t:main}, and let $\nu $ be as in \eqref{e:nudef}. Then, 
    for any $\eta >0$, there exists a constant,  
    $C=C(\beta ,\gamma,\nu,\eta)$, such that, $\mathcal Y$-a.s.~for 
    $N$ large enough, 
    for all $d\in\{0,\dots,N\}$
    \begin{equation}
      \label{e:rwa}
      \sum_{
        \substack{i,j=1\\
          \lfloor i/\nu \rfloor\neq 
          \lfloor j/\nu \rfloor}}
      ^{t r(N)}
      \bbone\{D_{ij}= d \}
      \le
      C\bigg[t^2 r(N)^2 2^{-N}\binom N d+
        t r(N) \nu^{-1}e^{\eta \|d\|}\bigg],
    \end{equation}
    \begin{equation}
      \label{e:rwc}
      \sum_{
        \substack{i,j=1, i\neq j\\
          \lfloor i/\nu \rfloor=
          \lfloor j/\nu \rfloor}}
      ^{t r(N)}
      \bbone\{D_{ij}= d \} 
      (\Lambda^0_d-\Lambda^1_{ij})
      \le 
      \frac{C d^2t r(N)}{ N^2} \bbone\{d\le \nu \}.
    \end{equation}
  \end{proposition}

  We now conclude the proof of Proposition~\ref{p:comp},
 that is we prove that the
  right-hand side of \eqref{e:comparison} tends to $0$. Observe first
  that $D_{ij}$ is smaller than $|i-j|$. Hence,  for 
  $\lfloor i/\nu \rfloor = \lfloor j/\nu \rfloor$
  \begin{equation}
    \Lambda^0_{ij}
    = \big[1-2N^{-1} D_{ij}\big]^p
    \ge [1-2N^{-1} |i-j|]^p \ge \Lambda^1_{ij}.
  \end{equation}
  Since $\Lambda_{ij}^1=0$ for  $(i,j)$ with  $\lfloor i/\nu \rfloor \neq
 \lfloor j/\nu \rfloor$, 
  $\Lambda_{ij}^0-\Lambda_{ij}^1< 0$ if and only if $\Lambda^0_{ij}<0$.  
  The summands on the right-hand side of \eqref{e:comparison} can be
  written as differences of two non-negative terms:
  \begin{equation}
    \label{e:twoterms}
    (\Lambda^0_{ij}-\Lambda^1_{ij})_+
    \mathbb E\Big[
      \frac{\partial^2 F_N(X^h_N)}{\partial X(i)\partial X(j)}
      \Big|\mathcal Y
      \Big]
    - 
    (\Lambda^0_{ij})_-
    \mathbb E\Big[
      \frac{\partial^2 F_N(X^h_N)}{\partial X(i)\partial X(j)}
      \Big|\mathcal Y
      \Big]
    .
  \end{equation}
  We bound this expression using Lemma~\ref{l:tech}. For given $\{u_i\}$
  let 
  \begin{equation}
    \tilde \Xi_N(c)=
    \max\{\bar \Xi_N(c,\beta ,\gamma ,u_i,u_j):1\le i,j\le \ell\}.
  \end{equation}
  Then $\tilde \Xi_N(c)$ satisfies \eqref{e:xibarprop} for some constants 
  $C$ and $C'$  and it is therefore increasing in $c$.   The absolute 
  value of the right-hand side of \eqref{e:comparison} is then bounded from 
  above by 
  \begin{equation}
    \begin{split}
      \label{e:bigest}
      \sum_{\substack{i,j =1\\i\neq j}}^{t r(N)}
      (\Lambda^0_{ij}&-\Lambda^1_{ij})_+
      \mathbb E\Big[
        \frac{\partial^2 F_N(X_N^0)}{\partial X(i)\partial X(j)}
        \Big | Y_N\Big]
      +
      \sum_{\substack{i,j =1\\i\neq j}}^{t r(N)}
      (\Lambda^0_{ij})_-
      \mathbb E\Big[
        \frac{\partial^2 F_N(X_N^1)}{\partial X(i)\partial X(j)}
        \Big]
      \\&\le
      \sum_{d=0}^N\Bigg\{
        \sum_{
          \substack{i,j=1\\
            \lfloor i/\nu \rfloor\neq 
            \lfloor j/\nu \rfloor}}
        ^{t r(N)}
        \bbone\{D_{ij}= d \} 
        (\Lambda^0_d)_+
        \int_{0}^1
        \tilde \Xi(h\Lambda^0_d)\d h
        \\&+
        \sum_{
          \substack{i,j=1,i\neq j\\
            \lfloor i/\nu \rfloor =
            \lfloor j/\nu \rfloor}}
        ^{t r(N)}
        \bbone\{D_{ij}= d \} 
        (\Lambda^0_d-\Lambda^1_{ij})
        \tilde \Xi\big(\Lambda^0_d\big)
        \\&+
        \sum_{i,j:|i-j|\ge N/2 }
        ^{t r(N)}
        \bbone\{D_{ij}= d \} 
        (\Lambda^0_d)_- 
        \tilde \Xi\big(0\big)
        \Bigg\}.
    \end{split}
  \end{equation}
  From the definition of $\tilde \Xi $ it follows that,
  \begin{equation}
    \int_0^1 \tilde \Xi(hc)\d h
    \le C\exp\Big\{-\frac{\gamma^2 N}{\beta^2 (1+c)}\Big\} 
    \int_0^1 (1-hc)^{-1/2}\d h.
  \end{equation}
  The last integral can be easily evaluated and is smaller than 2 for all 
  $c\in [-1,1]$.   
  Using Proposition~\ref{p:rwp}, the first line of \eqref{e:bigest}
  is smaller than the sum of the following two terms:
  \begin{equation}
    \label{e:terma}
    C\sum_{d=0}^N
    t^2 r(N)^2 2^{-N}\binom N d 
    \Lambda^0_d 
    \exp\Big\{-\frac{\gamma^2 N}{\beta^2 (1+\Lambda^0_d)}\Big\}
  \end{equation}
  and
  \begin{equation}
    \label{e:termb}
    C\sum_{d=0}^N
    \frac{t r(N)e^{\eta \|d\|}}{\nu }
    \Lambda^0_d 
    \exp\Big\{-\frac{\gamma^2 N}{\beta^2 (1+\Lambda^0_d)}\Big\}.
  \end{equation}
  The second line of \eqref{e:bigest} is bounded by 
  \begin{equation}
    \label{e:termc}
    C\sum_{d=0}^\nu 
    \frac{t r(N)d^2}{N^2}\,
    \tilde \Xi(\Lambda^0_d).
  \end{equation}
  The third line is non-zero only if $p$ is odd, and in that case it is 
  bounded by 
  \begin{equation}
    \label{e:termd}
    \sum_{d=N/2}^N C
    \bigg[t^2 r(N)^2 2^{-N}\binom N d+
      t r(N) \nu^{-1}e^{\eta \|d\|}\bigg]
    \Big(\frac {2d} N -1\Big)^p\tilde \Xi(0),
  \end{equation}

  We estimate \eqref{e:terma} first. Let
  $I(u)$ be  defined by
  \begin{equation}
    I(u)=u\log u + (1-u)\log (1-u) +\log 2,
  \end{equation}
  and let
  \begin{equation}
    J_N(u)=2^{-N}\binom N {\lfloor Nu\rfloor} \sqrt{\frac {\pi N}{2}}
    e^{NI(u)}.
  \end{equation}
  Stirling's formula yields 
  $J_N(u)\xrightarrow{N\to \infty} (4 u (1-u))^{-1}$ uniformly in $u$ on 
  compact subsets of $(0,1)$. Further, $J_N(u)\le C N^{1/2}$ for all 
  $u\in [0,1]$. From the definitions of $r(N)$ and $\tilde \Xi$, we find that 
  \begin{equation}
    \eqref{e:terma} =
    C\sum_{d=0}^N t^2 N^{1/2} 
    \Big(1-\frac{2d}N\Big)^p
    \exp\Big\{N \Upsilon_{p,\beta ,\gamma } \Big(\frac
        dN\Big)\Big\}J_N\Big(\frac dN\Big),
  \end{equation}
  where 
  \begin{equation}
    \Upsilon_{p,\beta ,\gamma }(u)=
    \begin{cases}
      \frac {\gamma^2}{\beta^2}-I(u)-
      \frac{\gamma^2}{\beta^2(1+(1-2u)^p)},&
      \text{if $(1-2u)^p\ge \frac \gamma {\beta^2}+\lambda-1 $,}
      \\
      \frac {\gamma^2}{\beta^2}-I(u)+
      \beta^2(1+(1-2u)^p)-2\gamma ,&
      \text{if $(1-2u)^p\le \frac \gamma {\beta^2}+\lambda-1 $.}
    \end{cases}
  \end{equation}

  \begin{lemma}
    \label{l:upsilon}
    There exists a function $\zeta (p)$ such that for all $p\ge 2$,  and 
    $\gamma $, $\beta $ satisfying 
    $\gamma \le \zeta (p) \beta $ and $\gamma < \beta^2$, 
    there exist positive constants  
    $\delta$, $ \delta '$ and 
    $c$ such that 
    \begin{equation}
      \label{e:upsa}
      \Upsilon_{p,\beta ,\gamma } (u)\le -\delta \qquad 
      \text{for all $u\in [0,1]\setminus(1/2-\delta' ,1/2+\delta' )$},
    \end{equation}
    and
    \begin{equation}
      \label{e:upsb}
      \Upsilon_{p,\beta ,\gamma } (u)\le -c (u-1/2)^2 \qquad
      \text{for all $u\in (1/2-\delta' ,1/2+\delta' )$}.
    \end{equation}
    Moreover $\zeta (p)$ is increasing and satisfies \eqref{e:zeta}, 
    that is
    \begin{equation}
      \label{e:zetab}
      \zeta (2)=2^{-1/2}, \quad 
      \zeta (3)=1.0291,
      \qquad \text{and} \qquad
      \lim_{p\to \infty}\zeta (p)=\sqrt{2\log 2}.
    \end{equation}
  \end{lemma}

  \begin{proof}
    Since $\gamma /\beta^2<1$, the second line of the definition of 
    $\Upsilon_{p,\beta ,\gamma } $ is used only for $p$ odd and 
    $u\ge u_c(p,\beta ,\gamma,\lambda  )=
    (1+(1-\lambda -\gamma /\beta^2)^{1/p})/2>1/2$.
    Furthermore,
 $\Upsilon_{p,\beta ,\gamma } (1/2)=\Upsilon_{p,\beta ,\gamma } '(1/2)=0$ and 
    \begin{equation}
      \Upsilon_{p,\beta ,\gamma } ''(1/2)=
      \begin{cases}
        4\big(\frac{2\gamma ^2}{\beta^2} -1\big),&
        \text{if $p=2$,}\\
        -4&
        \text{otherwise}.
      \end{cases}
    \end{equation}
    The second derivative is always negative for $\beta $, $\gamma $, $p$
    satisfying the assumptions of Theorem~\ref{t:main}. Therefore
    \eqref{e:upsb} holds.

    The second line of the definition of $\Upsilon_{p,\beta ,\gamma }(u)$ is
    decreasing in $u$. Hence for $u\ge u_c$
    \begin{equation}
      \Upsilon_{p,\beta ,\gamma }(u)\le 
      \Upsilon_{p,\beta ,\gamma }(u_c)=
      -\gamma (1-\gamma /\beta^2)-I(u_c)
    \end{equation}
    which is obviously strictly negative and \eqref{e:upsa} is proved for 
    $u\ge u_c$.  
    
    For any $\delta '>0$ and $u<1/2-\delta '$ the function $I(u)$ is
    strictly positive, and the function $\Phi (u)\equiv 
     1-1/(1+(1-2u)^p)$ is bounded. 
    Therefore, if $\gamma/\beta $ is sufficiently small, 
    then $\Upsilon_{p,\beta ,\gamma }(u)< - \delta $.
    If  $p$ is  even, 
    the function $\Upsilon_{p,\beta ,\gamma }$  is symmetric
    around $u=1/2$. 
    If $1/2<u<u_c(p,\beta ,\gamma )$ and $p$ is odd, then
    \begin{equation}
      \Upsilon_{p,\beta ,\gamma }(u)< \Upsilon_{p,1,0}(u)=-I(u)<0
    \end{equation}
    and the proof of \eqref{e:upsb} is finished.

    To prove the first part of \eqref{e:zetab} we should check that
    \eqref{e:upsb} holds for all $\gamma \le 2^{-1/2}\beta $.
    However, $\Upsilon_{2,\beta ,\gamma }(u)$
    is increasing in $\gamma^2/\beta^2$ and $I(u)\ge (1-2u)^2/2$.
    Thus, for $\gamma \le 2^{-1/2}\beta $,
    \begin{equation}
      \Upsilon_{2,\beta ,\gamma }(u)
      \le \frac 1 2 \Big(1-\frac 1{1+(1-2u)^2}\Big)-\frac 12(1-2u)^2.
    \end{equation}
    The right-hand side of the last inequality is equal $0$ for $u=1/2$ and
    its derivative
    \begin{equation}
      2(1-2u)\Big(1-\frac 1 {(1+(1-2u)^2)^2}\Big)>0 \quad 
      \text{for all $u<1/2$}.
    \end{equation}
    The symmetry of $\Upsilon_{2,\beta ,\gamma }$ around $1/2$ then implies
    the first part of \eqref{e:zetab}.

    Obviously, $\Phi (0)=1/2$, $\Phi '(0)=-2p$, $I(0)=\log2$ and 
    $I'(0)=-\infty$. Hence, for $\gamma/\beta =\sqrt{\log 2}$ there exists $u$
    small such that $\Upsilon_{p,\beta ,\gamma }(u)$ is positive. This
    implies $\zeta (p)<\sqrt{2\log 2}$.
    If $u\in (0,1/2)$ then $\lim_{p\to\infty}\Phi (u)=0$. This yield the
    second half of \eqref{e:zetab}.

    For illustration you find the graphs of function 
    $\Upsilon_{p,\beta ,\gamma }$ for $p=2,3,4$, $\beta =1$, and 
    $\gamma  =0$ 
    (solid lines), $\gamma =\sqrt{1/2}$ (dashed lines), $\gamma=1$
    (dash-dotted lines) and 
    $\gamma =\sqrt{2\log 2}$ (dotted lines) on  Figure~\ref{f:ups}. 
    The value of $\zeta (3)$ was calculated numerically using the figure
    for $p=3$.
    \begin{figure}
      \includegraphics[width=4.8cm]{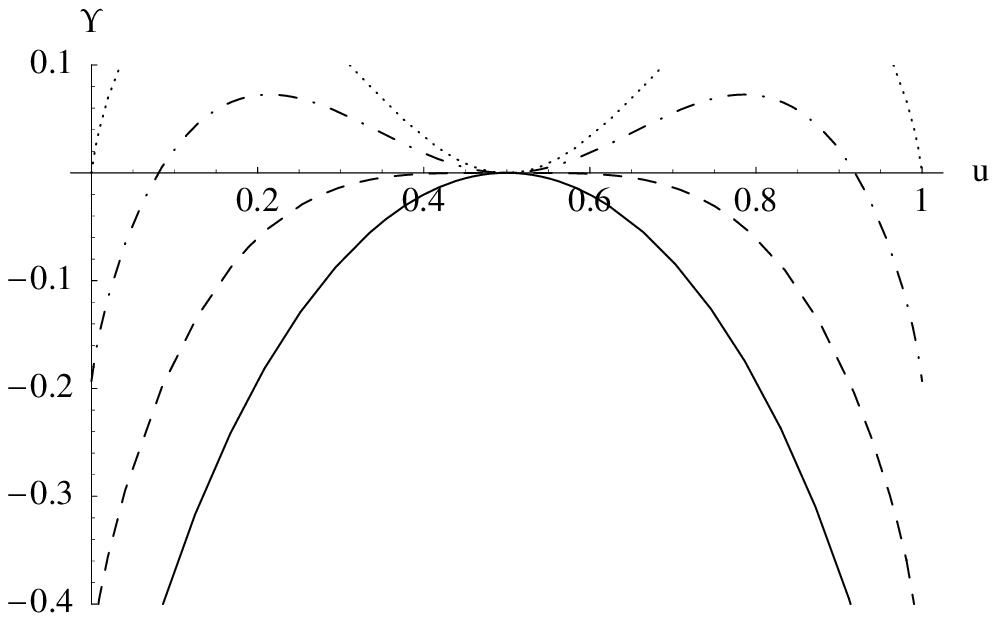}
      \includegraphics[width=4.8cm]{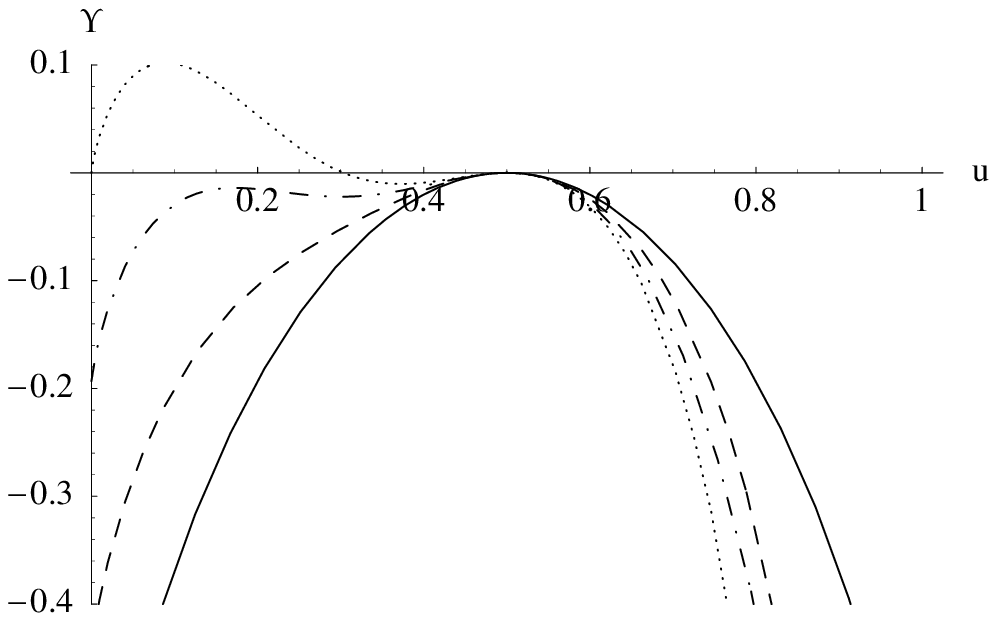}
      \includegraphics[width=4.8cm]{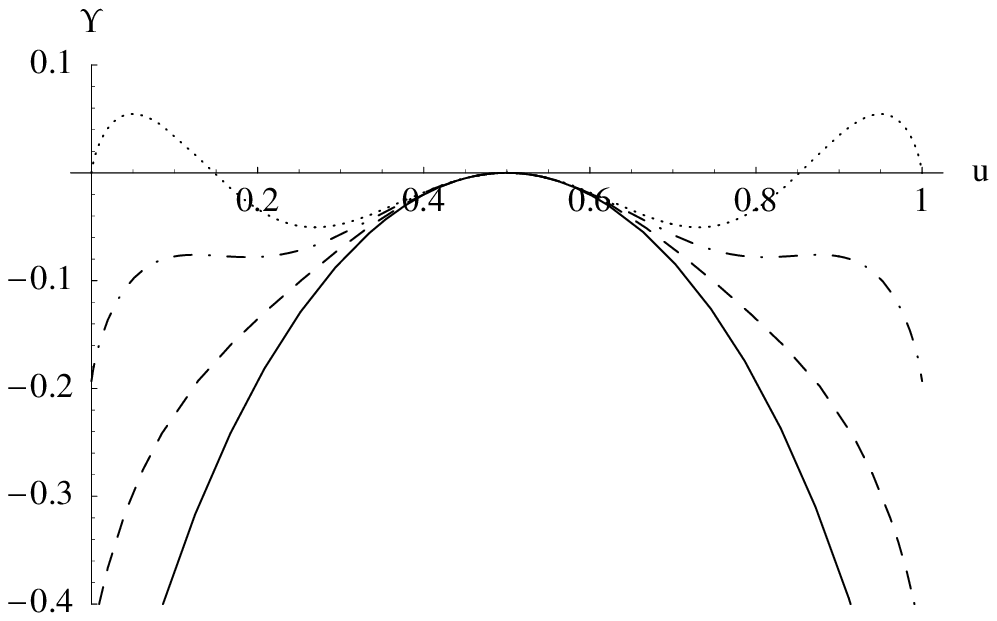}
      \caption{
        Function $\Upsilon_{p,\gamma ,\beta}$ for $p=2$, $3$, $4$ and
        various values of $\gamma/\beta$.}
      \label{f:ups}
    \end{figure}
  \end{proof}

  We can now finish the  bound  on \eqref{e:terma}. 
  Lemma~\ref{l:upsilon} and bounds on the function $J_N$ yield that
  for $d/N\notin (1/2-\delta' ,1/2+\delta' )$ the summands
  decrease exponentially in $N$. Therefore they can be neglected. 
  The remaining part can be bounded by 
  \begin{equation}
    \begin{split}
      C\sum_{d=(1/2-\delta ')N}^{(1/2+\delta ')N} 
      &t^2 N^{1/2}
      \Big(1-\frac{2d}N\Big)^p
      \exp(-cN (d/N-1/2)^2 )
      \\&\le
      C t^2 N^{3/2} \int_{-\delta '}^{\delta '}
      x^p e^{- c' N x^2}\d x 
      \\&\le
      C t^2 N^{3/2} N^{-(p+1)/2} 
      \int_{-\infty}^\infty u^p e^{-c'u^2}\d u
      \xrightarrow{N\to\infty} 0,
    \end{split}
  \end{equation}
  if $p\ge 3$.

  Similarly, for \eqref{e:termb} we have
  \begin{equation}
    \eqref{e:termb}\le
    C\sum_{d=0}^{N/2}
    t N^{1/2}\nu^{-1} 
    \Big(1-\frac{2d}N\Big)^p
    \exp(N\tilde\Upsilon(d/N)),
  \end{equation}
  where, setting $\|u\|=\min (u,1-u)$,
  \begin{equation}
    \tilde\Upsilon_{p,\beta ,\gamma }(u)=
    \begin{cases}
      \frac {\gamma^2}{2\beta^2}-
      \frac{\gamma^2}{\beta^2(1+(1-2u)^p)}+\eta \|u\|,&
      \text{if $(1-2u)^p\ge \frac \gamma {\beta^2}+\lambda -1$,}
      \\
      \frac {\gamma^2}{2\beta^2}+
      \beta^2(1+(1-2u)^p)-2\gamma+\eta \|u\|,&
      \text{if $(1-2u)^p\le \frac \gamma {\beta^2}+\lambda -1$.}
    \end{cases}
  \end{equation}
  Observe first that the second part of the definition of 
  $\tilde \Upsilon_{p,\beta ,\gamma }$ is always strictly negative.
  It is also easy to be checked that it is possible to choose 
  $\delta$, $\delta '$ and $\eta $ small such
  that the first part of the definition of 
  $\tilde \Upsilon (u)<\delta $ for all $\|u\|\ge \delta '$.
  Therefore such $d$ can be neglected. Around $d=0$ 
  the function $\tilde\Upsilon (x)$ can be approximated by a 
  linear function $-c x$, $c>0$, and the summation by an integration.
  As an upper bound we get
  \begin{equation}
    \label{e:termbint}
    Ct N^{3/2} \nu^{-1} \int_0^{\delta '}
    e^{-c N x} \d x \le C t N^{1/2}\nu^{-1}
    \xrightarrow{N\to\infty} 0.
  \end{equation}
  An analogous bound works for $d$ close to $N$ and $p$ even.

  For \eqref{e:termc} we have
  \begin{equation}
    \label{e:termcsum}
    \eqref{e:termc} \le 
    C\sum_{d=0}^\nu 
    t N^{-3/2} d^2 
    [1-(1-2dN^{-1})^p]^{-1/2}
    \exp(N\tilde\Upsilon (d/N)).
  \end{equation}
  The linear approximation of $\tilde\Upsilon $ and of the bracket in the
  last expression  yields an upper bound
  \begin{equation}
    \label{e:termcint}
    C t N^{3/2}  \int_0^\varepsilon x^{3/2} e^{-c'N x}\d x 
    \le 
    C t N^{-1}
    \xrightarrow{N\to\infty} 0.
  \end{equation}

  Finally, since $\tilde \Xi (0)=C e^{-N\gamma^2/\beta^2}$, 
  it is easy to see that the second half of \eqref{e:termd} tends to $0$. 
  The first half equals (up to constant) 
  \begin{equation}
    \begin{split}
      \label{e:kla}
      \sum_{d=N/2}^N&
      \Big(\frac{2d}N-1\Big)^p t^2 N 2^{-N}\binom N d
      \\&\le 
      Ct^2\Big\{\sum_{d\ge N/2 + N^{3/5}}  N 2^{-N}\binom N d
        +\sum_{i=1}^{2N^{3/5}}\Big(\frac {N+i} N -1\Big)^p N^{1/2} e^{-i^2/2N}
        \Big\},
    \end{split}
  \end{equation}
  where we used the known approximation of 
  $\binom N d \le C N^{-1/2}2^N e^{-i^2/2N}$ for $d=(N+i)/2$ and 
  $i\ll N^{2/3}$. The first term in \eqref{e:kla} tends to $0$ by a
  standard moderate deviation argument. The second one can be approximated
  by 
  \begin{equation}
    Ct^2 N^{1-(p/2)}\int_{0}^\infty x^p e^{-x^2/2}\d
    x\xrightarrow{N\to\infty}0
  \end{equation}
  for $p\ge 3$.
  This completes the proof of Proposition~\ref{p:comp}.
\end{proof}

\section{Random walk properties}
\label{s:rwprop}
In this section we prove Proposition~\ref{p:rwp}.
For $A\subset \mathcal S_N$ let 
$T_A=\min\{k\ge 1: Y_N(k)\in A\}$ be the hitting time of $A$. We write 
$\mathbb P_x$ for the law of the simple random walk $Y_N$ conditioned on 
$Y_N(0)=x$. 
Let $Q=Q_i$, $i\in \mathbb N$, be a birth-death process on $\{0,\dots,N\}$
with transition probabilities $p_{i,i-1}=1-p_{i,i+1}=i/N$. 
We use 
$P_k$ and $E_k$ to denote the law of (the expectation with respect to)
$Q$ conditioned on $Q_0=k$.
Under $P_0$,
$Q_i$ has the same law as $\dist(Y_N(0),Y_N(i))$.
Define $T_k=\min\{i\ge 1: Q_i=k\}$ the
hitting time of $k$ by $Q$.
It is well-known fact that for $k<l<m$
\begin{equation}
  \label{e:ehr}
  P_l[T_m<T_k]=
  \frac{\sum_{i=k}^{l-1}\binom{N-1}i^{-1}}
  {\sum_{i=k}^{m-1}\binom{N-1}i^{-1}} .
\end{equation}
Finally, let $p_k(d)=P_0(Q_k=d)$.  
We need the following lemma for estimating $p_k(d)$ for large $k$. 
\begin{lemma}
  \label{l:distance}
  There exists $K$ large enough such that  for all 
  $k\ge K N^2 \log N=:\mathcal K(N)$ and 
  $x,y\in \mathcal S_N$ 
  \begin{equation}
    \bigg| \frac{\mathbb P_y[Y_N(k)=x\cup Y_N({k+1})=x]}2-2^{-N}\bigg|
    \le 2^{-8 N}
  \end{equation}
  and thus
  \begin{equation}
    \bigg| \frac{p_k(d)+p_{k+1}(d)}2-2^{-N}\binom N d\bigg|\le 2^{-4 N}.
  \end{equation}
\end{lemma}

\begin{proof}
  The beginning of the argument is the same as in \cite{Mat87}. We 
  construct coupling between $Y_N$ (which by definition starts at site 
    $\boldsymbol 1=(1,\dots,1)\in \mathcal S_N$)  
  and another process $Y_N^\star$. This
  process is a simple random walk on $\mathcal  S_N$ 
  with the initial distribution $\mu_N^\star$ being uniform 
  on those $x\in \mathcal S_N$ with 
  $\dist(x,\boldsymbol 1)$ even. The coupling is the same as  in \cite{Mat87}.
  This coupling gives certain random time $\mathcal T_N$ which can be used to bound
  the variational distance between $\mu^{\star}$ and the distribution 
  $\mu^k_N$ of 
  $Y_N(k)$: for $k$ even
  \begin{equation}
    d_\infty(\mu^\star_N,\mu^k_N)
    \equiv
    \max_{A\subset \mathcal S_N}
    |\mu^\star_N(A)-\mu^k_N(A)|
    \le \mathbb P[\mathcal T_N> k].
  \end{equation}
  The law of $\mathcal T_N$ is as follows. Let 
  $U=\dist(Y_N^\star(0),\boldsymbol 1)$. That is $U$ is a binomial random 
  variable with parameters $N$ and $1/2$ conditioned on being even. 
  Consider another simple random walk $\tilde Y_U$ on $\mathcal S_U$ 
  started from $\boldsymbol 1$. The distribution of $\mathcal T_N$ is then 
  the same as the distribution of the hitting time of 
  $\{x\in \mathcal S_U:\dist(\boldsymbol 1,x)=U/2\}$. It is proved in 
  \cite{Mat87} that $P(\mathcal T_N>N\log N)\to c < 1$. It is then easy to 
  see that,
  \begin{equation}
    \mathbb P[\mathcal T_N\ge \mathcal K(N)]\le c^{KN/2}\le 2^{-8N},
  \end{equation}
  if $K$ is large enough. Thus, for even $k\ge \mathcal K(N) $, 
  $d_{\infty}(\mu_N^\star,\mu_N^k)\le 2^{-8N}$ and thus 
  $|\mu_N^\star(x)-\mu_N^k(x)|\le 2^{-8N}$ for all $x\in \mathcal S_N$. A 
  similar claim for $k$ odd is then not difficult to prove. The second 
  part of the lemma is a direct consequence of the first part. 
\end{proof}

\begin{lemma}
  \label{l:rwpp}
  Let $\gamma $, $\beta $, $\nu $ satisfy the hypothesis of 
  Proposition~\ref{p:rwp}. Then, 
  there exists a constant,  
  $C=C(\beta ,\gamma,\nu)$, such that  for all $N$ large enough, 
  $\mathcal Y$-a.s.
  \begin{equation}
    \label{e:rwb}
    \sum_{
      \substack{i,j=1,i\neq j\\
        \lfloor i/\nu \rfloor=
        \lfloor j/\nu \rfloor}}
    ^{t r(N)}
    \bbone\{D_{ij}= d \}
    \le 
    C t r(N)\bbone\{d\le \nu \},
  \end{equation}
  and for all $d\in\{0,\dots,N\}$.
\end{lemma}
\begin{proof}
  The lemma is trivially true for $d>\nu $. For $d\le \nu $, let
  \begin{equation}
    \rho (d)=E_0\sum_{i=1}^\nu \bbone\{Q_i=d\}.
  \end{equation}
  We have $\rho (0)\ge N^{-1}$ and 
  $\rho (d)\ge P_0[T_d\le \nu ]$. This probability is decreasing in 
  $d$ and 
  \begin{equation}
    \label{e:rovne}
    P_0[T_\nu \le \nu ]=\frac NN \cdot \frac {N-1}N \dots
    \frac{N-\nu +1}{N}
    \ge e^{-\nu^2/N}.
  \end{equation}
  Thus $\rho (d)\ge e^{-\nu^2/N}$ for all $d\le \nu $.
  To get an upper bound on $\rho (d)$ we write
  \begin{equation}
    \rho (d)\le E_0\Big[\sum_{i=1}^{T_\nu }\bbone\{Q_i=d\}\Big]
    =1+E_d\Big[\sum_{i=1}^{T_\nu }\bbone\{Q_i=d\}\Big]=
    1+\frac 1 {P_d[T_\nu < T_d]}.
  \end{equation}
  However, using \eqref{e:ehr},
  \begin{equation}
    P_d[T_\nu < T_d]
    =\frac{N-d}{N}P_{d+1}[T_\nu <T_d]
    =\frac{N-d}{N}
    \frac{\binom{N-1}{d}^{-1}}
    {\sum_{i=d}^{\nu -1}\binom{N-1}i^{-1}}
    =1-O(\nu N^{-1} ).
  \end{equation}
  Since $\nu \ll N$, $\rho (d)\le 2$.

  Consider now one-block contribution to \eqref{e:rwb}, 
  \begin{equation}
    \sum_{i,j=1}^\nu \bbone\{D_{ij}=d\} =:\nu^2 \tilde Z. 
  \end{equation}
  Of
  course, $\tilde Z\in [0,1]$ and, using the results of the previous
  paragraph,
  \begin{equation}
    \label{e:EZ}
    e^{-\nu^2/N} (2\nu )^{-2}  \le  \mathbb E[\tilde Z]\le 2\nu^{-1}. 
  \end{equation}
  The left-hand side of \eqref{e:rwb} is stochastically smaller than
  $ \nu^2 \sum_{k=1}^{m}\tilde Z_k$,
  where $\tilde Z_k$ are i.i.d.~copies of $\tilde Z$ and 
  $m=\lceil tr(N)/\nu \rceil$. By Hoeffding's
  inequality \cite{Hoe63},
  \begin{equation}
    \mathbb P\Big[\sum_{i=1}^{m} \tilde Z_k 
      \ge 2 m  \mathbb E[\tilde Z_k]\Big]
    \le \exp\{-2 m^2  \mathbb E[\tilde Z_k]^2\} 
    \le \exp\{-m^2 e^{-2\nu^2/N} (2\nu )^{-4}\},
  \end{equation}
  where we used the lower bound from \eqref{e:EZ}. Since $\nu /N^2\ll N$,
  by the
  Borel-Cantelli lemma, the left-hand side
  of \eqref{e:rwb} is a.s.~bounded by
  \begin{equation}
    \nu^2 2 m \mathbb E[\tilde Z]\le C t r(N)
  \end{equation}
  for all $N$ large enough and $d\le \nu $. This completes the
  proof of Lemma~\ref{l:rwpp}.
\end{proof}

\begin{proof}[Proof of Proposition~\ref{p:rwp}]
  We prove \eqref{e:rwc} first.
  Observe that for $i,j$ in the same block
  \begin{equation}
    \label{e:rwhb}
    \Lambda^0_d - \Lambda^1_{ij}=
    \Big(1-\frac{2d}N\Big)^p
    -\Big(1-\frac{2p|i-j|}{N}\Big)=
    \frac{2p(|i-j|-d)}{N}+O\Big(\frac {d^2}{N^{2}}\Big).
  \end{equation}
  The contribution of the error term is smaller
  than the right-hand side of \eqref{e:rwc}, as follows from
  Lemma~\ref{l:rwpp}. 

  To compute the contribution of the main term, let
  \begin{equation}
    \tilde \rho (d)=E_0\Big[\sum_{i=1}^\nu 
      (i-d) \bbone\{Q_i=d\}\Big]. 
  \end{equation}
  Let $T_d^{1}=T_d$ and $T_d^k=\min\{i> T_d^{k-1}:Q_i=d\}$. Then
  \begin{equation}
    \begin{split}
      \label{e:trho}
      \tilde \rho (d)
      &=E_0 \Big[
        \sum_{j=1}^\infty( T_d^j-d)\bbone\{T_d^j<\nu \}\Big]
      =E_0 \Big[
        \sum_{j=1}^\infty( T_d^j-T_d^1+T_d^1-d)\bbone\{T_d^j<\nu \}\Big]
      \\&\le 
      E_0[(T_d -d) \bbone\{T_d<\nu \}]\Big(1+
        \sum_{i=1}^\infty E_d[T_d^i\bbone\{T_d^i<\nu-d \}]\Big).
    \end{split}
  \end{equation}
  Using \eqref{e:rovne}, $P_0[T_d=d]\le C e^{-d^2/N}$ and further
  \begin{equation}
    P_0[T_d\ge d+2k]\le \binom{d+2k}{k}
    \Big(\frac d N\Big)^k\le C\frac {d^{2k}}{N^k}.
  \end{equation}
  Hence, $ cd^2 N^{-1}\le E_0[(T_d -d) \bbone\{T_d<\nu \}]\le C  d^2 N^{-1}$.

  For the second term in \eqref{e:trho} we write
  \begin{equation}
    \begin{split}
      \label{e:rwhh}
      1+&
        \sum_{i=1}^\infty E_d[T_d^i\bbone\{T_d^i<\nu-d \}]
      \\&\le
      1+E_d[T_d\bbone\{T_d<\nu-d \}]\Big(1+
          \sum_{i=1}^\infty E_d[T_d^i\bbone\{T_d^i<\nu-d \}]\Big)
      \\&=
      \sum_{k=0}^\infty
      \big\{E_d[T_d\bbone\{T_d<\nu-d \}]\big\}^k.
    \end{split}
  \end{equation}
  Using the well-known estimate $\binom {2k} k \le C k^{-1/2} 2^k$ and $k<2^k$,
  \begin{equation}
    E_d[T_d\bbone\{T_d<\nu-d \}]\le 
    \sum_{k=1}^{\nu /2}2k \binom {2k} k \Big(\frac \nu N\Big)^k\le
    C\sum_{k=1}^{\infty}\Big(\frac {4 \nu }N\Big)^k\le
    C \frac \nu  N
  \end{equation}
  and \eqref{e:rwhh} is finite. Thus $\tilde \rho (d)\le C d^2 N^{-1}$ for all 
  $d\in \{0,\dots , \nu \}$.

  The one-block contribution of the first term of~\eqref{e:rwhb} to 
  \eqref{e:rwc} is then given by 
  \begin{equation}
    \frac {2p}{N}\sum_{i,j=1}^\nu (|i-j|-d)\bbone \{D_{ij}=d\}=:\frac{2p}N\nu^3\tilde Z,
  \end{equation}
  with $\tilde Z\in [0,1]$ and 
  \begin{equation}
    \label{e:EZb}
    c d^2 N^{-1}\nu^{-3} \le \mathbb E[\tilde Z] \le C d^2 N^{-1} \nu^{-2}.  
  \end{equation}
  Therefore, as in the proof of Lemma~\ref{l:rwpp},  Hoeffding's 
  inequality and \eqref{e:EZb}  imply that 
  the contribution of the first term of~\eqref{e:rwhb} to \eqref{e:rwc}
  is smaller than $C t r(N) d^2 N^{-2}$, which was to be shown.

  Finally, we prove \eqref{e:rwa}. Since we are interested in an upper
  bound only we can,  
  without loss of generality, restrict the summation on $i<j$. We
  first consider the contribution of pairs $(i,j)$ such that 
  $j-i\ge \mathcal K(N)$. Then necessarily, 
  $\lfloor i/n\rfloor\neq\lfloor j/n\rfloor$. 
  Let $R=tr(n)$.
  Lemma~\ref{l:distance} yields 
  \begin{equation}
    \mathbb E\Big[\sum_{j-i\ge \mathcal K(N)}^R \bbone\{D_{ij}=d\}\Big]
    = \sum_{j-i\ge \mathcal K(N)}^R p_{j-i}(d)\le CR^2 2^{-N}\binom N d.
  \end{equation}
  Further,
  \begin{equation}
    \begin{split}
      \label{e:varest}
      &\Var \bigg[
        \sum_{j-i\ge \mathcal K(N)}^{R}
        \bbone\{D_{ij}= d \}\bigg]
      \\&=
      \sum_{j_1-i_1\ge \mathcal K(N)}^R
      \sum_{j_2-i_2\ge \mathcal K(N)}^R
      \mathbb
      P\big[D_{i_1,j_1}=D_{i_2,j_2}=d\big]
      -
      \mathbb P\big[D_{i_1,j_1}=d\big]
      \mathbb P\big[D_{i_2,j_2}=d\big].
    \end{split}
  \end{equation}
  We can again suppose that $i_1\le i_2$. The right-hand side of
  \eqref{e:varest} is non-null only if  
  $i_1\le i_2 \le j_1<j_2$ or $i_1\le i_2 < j_2 \le j_1$. We will
  consider only the first case. The second one can be treated
  analogously. In is not difficult to see using Lemma~\ref{l:distance}
  that if $i_2-i_j\ge \mathcal K(N)$ or $j_2-j_1\ge \mathcal K(N)$ then the difference of
  probabilities in the above summation is at most $2^{-4N}$. Therefore,
  the contribution of such $(i_1,i_2,j_1,j_2)$ to the variance is at most
  $R^4 2^{-4N}$. 

  If $i_2-i_1<\mathcal K(N)$ and $j_2-j_1<\mathcal K(N)$ then, using Lemma~\ref{l:distance}
  again,  
  \begin{equation}
    \mathbb P\big[D_{i_1,j_1}=D_{i_2,j_2}=d\big]
    \le C2^{-N}\binom N d.
  \end{equation}
  We choose $\varepsilon >0$. For $\|d\|\le (1-\varepsilon )N/2$ we have
  \begin{equation}
    \begin{split}
      \sum_{\substack{
          j_1-i_1\ge \mathcal K(N)\\
          i_2-i_1< \mathcal K(N)
      }}
      &
      \sum_{\substack{
          j_2-i_2\ge \mathcal K(N)\\
          j_2-j_1< \mathcal K(N)
      }}
      \mathbb
      P\big[D_{i_1,j_1}=D_{i_2,j_2}=d\big]
      \\&
      \le C \mathcal K(N)^2 R^2 2^{-N}\binom N d 
      \le C \mathcal K(N)^2 R^2 e^{-NI((1-\varepsilon /2)/2)} \ll N^{-3} R^2 \nu^{-2},
    \end{split}
  \end{equation}
  say. For $\|d\|\ge (1-\varepsilon )N/2$, that is 
  $|d-N/2|\le \varepsilon N/2$, we have for $\varepsilon $ small enough
  (how small depend on $\gamma$ and $\beta $) that 
  $2^{-N}\binom N d \gg N^7 R^{-2}$. Then,
  \begin{equation}
    \begin{split}
      \sum_{\substack{
          j_1-i_1\ge \mathcal K(N)\\
          i_2-i_1< \mathcal K(N)
      }}
      &
      \sum_{\substack{
          j_2-i_2\ge \mathcal K(N)\\
          j_2-j_1< \mathcal K(N)
      }}
      \mathbb
      P\big[D_{i_1,j_1}=D_{i_2,j_2}=d\big]
      \\&
      \le 
      C N^4 R^2  2^{-N}\binom Nd \ll N^{-3}R^4 2^{-2N}\binom N d^2.
    \end{split}
  \end{equation}
  We have thus found that the expectation of the summation over 
  $j-i>\mathcal K(N)$ is smaller than the right-hand side of \eqref{e:rwa} and the
  variance of the same summation is much smaller than $N^{-3}$ times the
  right-hand side of \eqref{e:rwa} squared. A straightforward application
  of the Chebyshev inequality and  the Borel-Cantelli Lemma then gives
  the desired a.s.~bound for pairs $j-i\ge \mathcal K(N)$ and all 
  $d \in \{0,\dots,N\}$.

  Choose again $\varepsilon >0$. For $j-i<\mathcal K(N)$, observe first that if 
  $\|d\|\ge (\log N)^{1+\varepsilon }\gg \log N$ then the summation 
  over such pairs $(i,j)$ in \eqref{e:rwa} is always smaller than 
  $\mathcal K(N) R \ll R\nu^{-1}e^{\eta \|d\|}$ for all $\eta >0$. For the remaining 
  $d$'s, that is $\|d\|<(\log N)^{1+\varepsilon '}$, let $K_N\ge K$ be 
  the smallest constant such that $K_N N^2 \log N$ is a multiple of $\nu $. 
  Since $\nu \ll N^2$, $K_N-K\ll 1$. As the difference between $K$ 
  and $K_N$ is negligible, we will use the same notation 
  $\mathcal K(N)$ for $K_N N^2\log N$ and we will  
  simply suppose that $\mathcal K(N)$ is a multiple of $\nu $. The summation in 
  \eqref{e:rwa} for $j-i\le \mathcal K(N)$ can be bounded from above by 
  \begin{equation}
    \sum_{
      \substack{0<j-i<\mathcal K(N)\\
        \lfloor i/\nu \rfloor\neq 
        \lfloor j/\nu \rfloor}}
    ^{t r(N)}
    \bbone\{D_{ij}= d \}
    \le 
    \sum_{k=0}^{\mathcal K(N)-1}\sum_{\ell =0}^{\lceil R/\mathcal K(N)\rceil}
    \sum_{m=j_{k}}^{\mathcal K(N)}
    \bbone\{D_{\mathcal K(N)\ell+k,\mathcal K(N)\ell+k+m}= d \},
  \end{equation}
  where $j_{k}$ is the smallest integer such that 
  $\lfloor (\mathcal K(N)\ell+k)/\nu\rfloor\neq \lfloor (\mathcal K(N)\ell+k+j_{k})/\nu\rfloor$, 
  which does not depend on $\ell$.
  We define random variables $Z_\ell(j,d)$ by
  \begin{equation}
    Z_\ell(j,d)=\frac 1 {\mathcal K(N)}
    \sum_{m=j}^{\mathcal K(N)}
    \bbone\{D_{\mathcal K(N)\ell+k,\mathcal K(N)\ell+k+m}= d \}.
  \end{equation}
  The sequence $\{Z_\ell(j,d):\ell\ge 0\}$ for fixed $j$ and $d$ is a sequence
  of i.i.d.~variables with values in $[0,1]$.

  Let 
  $E_N=\{d:\|d\|<(\log N)^{1+\varepsilon '}, d\ge N/2\}$. For $d\in E_N$
  \begin{equation}
    \mathbb P[Z_\ell(k,d)>0]
    \le \binom N {d} P_{d}(T_{\boldsymbol 1}<\mathcal K(N))
    \le \binom N {d} e^{\lambda K} 
     E_{d}\big[e^{-\lambda T_{\boldsymbol 1}/N^2\log N }\big].
  \end{equation}
  According to Lemma 3.4 of \cite{CG06},
  \begin{equation}
    E_{d}\big[\exp(-\lambda T_{\boldsymbol 1}m(N)^{-1})\big]
    \le (2^{-N} m(N) \lambda^{-1}+\xi_N(d))(1+o(1)),
  \end{equation}
  for $N\log N \ll m(N)\ll 2^N$, with 
  $\xi_n(k)=2^{-n } \frac n 2 \binom n k^{-1} \sum_{j=1}^{n-k}\binom n {k+j} \frac 1 j$.
  Taking $m(N)=N^2$ and $d\in E_N$ it is not difficult to check that for 
  $\varepsilon $ small enough
  \begin{equation}
    \mathbb E_{z_d}\big[e^{-\lambda T_{\boldsymbol 1}/N^2 }\big]\le
    2^{-N(1-\varepsilon )}.
  \end{equation}
  Hence,
  \begin{equation}
    \begin{split}
      \mathbb P\Big[
        \bigcup_{d\in E_N} \Big\{
          \sum_{k=0}^{\mathcal K(N)-1}\sum_{\ell =0}^{\lceil R/\mathcal
            K(N)\rceil}
          &
          Z_\ell(j_k,d)>0\Big\}\Big]
      \\&
      \le C 
      \binom N {\lceil (\log N)^{1+\varepsilon }\rceil} 
      R (\log N)^{1+\varepsilon }2^{-N(1-\varepsilon )}
      \le C 2^{-\varepsilon ' N },
    \end{split}
  \end{equation}
  for some $\varepsilon '$ small.
  Hence, $d\in E_N$ do not pose any problem, by the
  Borel-Cantelli lemma again.

  To treat $d\le (\log N)^{1+\varepsilon' }$  we will distinguish two
  cases: $j_k\le 2d$ and $j_k>2d$. For the first case,
  observe that for any $d<\nu $ there
  are at most $d \mathcal K(N)/\nu $ values of 
  $k\in \{0,\dots, \mathcal K(N)-1\}$ 
  such that $j_k\le d$. 
  Further, as before, $Z_{\ell}(j_k,d)\le  Z_\ell(0,d)$, 
  $\mathbb E[ Z_\ell(0,d)]\ge 1/(N\mathcal K(N))$, and 
  $\mathbb E[Z_\ell(0,d)]\le C/\mathcal K(N)$. Hence, by Hoeffding's
  inequality, the probability
  \begin{equation}
    \mathbb P\Big[\mathcal K(N) \sum_{\ell =0}^{\lceil R/\mathcal K(N)\rceil}
      Z_\ell(0,d) \ge \frac R{\mathcal K(N)}\Big]
  \end{equation}
  decreases at least exponentially with $N$ and thus for $j_k< 2d$, a.s,
  \begin{equation}
    \label{e:poa}
    \mathcal K(N) \sum_{\ell =0}^{\lceil R/\mathcal K(N)\rceil}
    Z_\ell(0,d) \ge \frac R{\mathcal K(N)}.
  \end{equation}
  For $j\ge 2d$ and $N$ large enough, $Z_{\ell}(j,d)\le Z_\ell(d+6,d)$. 
  We have,
  \begin{equation}
    c N^{-6}\le \mathcal K(N) \mathbb E[Z_\ell(d+6,d)] \le C N^{-3}.
  \end{equation}
  Indeed, the lower bound is trivial and for the upper bound we use the 
  fact that the probability that $Y_N$ reaches $d$ before returning to 
  $d+6$ is smaller than $CN^{-5}$ and before the time $\mathcal K(N)$ there are at 
  most $\mathcal K(N)$ tries. Hence, for $j\ge 2d$ the probability
  \begin{equation}
    \mathbb P\Big[\mathcal K(N) \sum_{\ell =0}^{\lceil R/\mathcal K(N)\rceil}
      \tilde Z_\ell(k,d) \ge \frac R{N^3\mathcal K(N)}\Big]
  \end{equation}
  decreases at least exponentially in $N$ and thus the interior inequality is
  not valid a.s.~for all $N$ large.
  Summing over $k$ we 
  get 
  \begin{equation}
    \sum_{k=0}^{\mathcal K(N)-1}\sum_{\ell =0}^{\lceil R/\mathcal
      K(N)\rceil} \mathcal K(N)
    Z_\ell(j_k,d)\le d \mathcal K(N) \nu^{-1} \frac R{\mathcal K(N)}+
    \mathcal K(N) \frac R{N^3\mathcal K(N)} \le C R\nu^{-1}e^{\eta d},
  \end{equation}
  since $\gamma /\beta^2 <1$.
\end{proof}

\section{Convergence of clock process} 
\label{s:clock}

We will prove the convergence of the rescaled clock process  to the stable
subordinator on space $D([0,T],\mathbb R)$ equipped with the Skorokhod
$M_1$-topology. This topology is not commonly  used in the literature,
therefore we shortly recall some of its properties and compare it with the more
standard Skorokhod $J_1$-topology, which we will need later, too. For more
details the reader is referred  to \cite{Whi02} for both topologies and to 
\cite{Bil68} for detailed account on $J_1$-topology. 

\subsection{Topologies on the Skorokhod space }
Consider space $D=D([0,T],\mathbb R)$ of càdlàg functions. The $J_1$-topology
is the topology given by the $J_1$-metric: for $f,g\in D$
\begin{equation}
  d_{J_1}(f,g)=\inf_{\lambda \in \Lambda }
  \{\|f\circ \lambda -g\|_\infty \vee \|\lambda -e\|_\infty\},
\end{equation}
where $\Lambda $ is the set of strictly increasing functions mapping $[0,T]$
onto itself such that both $\lambda $ and its inverse are continuous, and 
$e$ is the identity map on $[0,T]$.  

Also the $M_1$-topology is given by a metric. For $f\in D$ let 
$\Gamma_f$ be its completed graph,
\begin{equation}
  \Gamma_f=\{(z,t)\in \mathbb R\times[0,T]: z=\alpha f(t-)+(1-\alpha)f(t),
    \alpha \in [0,1]\}.
\end{equation}
A parametric representation of the completed graph $\Gamma_f$ (or of $f$)
is a continuous bijective mapping $\phi(s)=(\phi_1(s),\phi_2(s))$, 
$[0,1]\mapsto\Gamma_f$ whose first   coordinate 
$\phi_1$ is increasing. If 
$\Pi (f)$ is set of all parametric representation of $f$, then the
$M_1$-metric is defined by
\begin{equation}
  d_{M_1}(f,g)=\inf\{\|\phi_1-\psi_1\|_\infty
    \vee \|\phi_2-\psi_2\|_\infty: \phi \in \Pi (f), \psi \in \Pi (g)\}.
\end{equation}
The space $D$ equipped with both $M_1$- and $J_1$-topologies is Polish. The
$M_1$-topology is weaker than the $J_1$-topology: As an example, consider 
the sequence 
\begin{equation}
  f_n=\bbone\{[1-1/n,1)\}+ 2\cdot \bbone\{[1,T]\}, 
\end{equation}
which converges to 
$f=2\cdot \bbone\{[1,T]\}$ in the $M_1$-topology but not in the $J_1$-topology. One
often says that the $M_1$-topology allows ``intermediate jumps''. 

We will need a criterion for tightness of probability measures on $D$. To this
end we define several moduli of continuity,
\begin{equation}
  \begin{split}
    w_f(\delta )&
    =\sup\big\{ \min\big(|f(t)-f(t_1)|,|f(t_2)-f(t)|\big): 
      t_1\le t\le t_2\le T,t_2-t_1\le \delta \big\},
    \\
    w'_f(\delta )&=
    \sup \big\{ \inf_{\alpha \in [0,1]}|f(t)-(\alpha f(t_1)+(1-\alpha )f(t_2))|:
      t_1\le t\le t_2\le T,t_2-t_1\le \delta \big\},
    \\
    v_f(t,\delta )&=
    \sup\big\{|f(t_1)-f(t_2)|:t_1,t_2\in [0,T]\cup (t-\delta ,t+\delta )\big\}.
  \end{split}
\end{equation}
The following result is a restatement of Theorem 12.12.3 of \cite{Whi02}
and Theorem~15.3 of \cite{Bil68}.
\begin{theorem}
  \label{t:tight}
  The sequence of probability measures $\{P_n\}$ on $D([0,T],\mathbb R)$ 
  is tight in the $J_1$-topology if 
  \begin{enumerate}
    \item[(i)] For each positive $\varepsilon $ there exist $c$ such that
    \begin{equation}
      P_n [f:\|f\|_\infty>c]\le \varepsilon , \qquad n\ge 1.
    \end{equation}
    \item[(ii)] For each $\varepsilon >0$ and $\eta >0$, there exist 
    a $\delta $, $0<\delta <T$, and an integer $n_0$ such that
    \begin{equation}
      \label{e:Jcond}
      P_n[f:w_f (\delta )\ge \eta ]\le \varepsilon , \qquad n\ge n_0,
    \end{equation}
    and
    \begin{equation}
      P_n[f:v_f (0,\delta )\ge \eta ]\le \varepsilon \text{ and }
      P_n[f:v_f (T,\delta )\ge \eta ]\le \varepsilon , \qquad n\ge n_0.
    \end{equation}
  \end{enumerate}

  The same claim hold for the $M_1$-topology with $w_f(\delta )$ in
  \eqref{e:Jcond} replaced by $w'_f(\delta )$.
\end{theorem}

\subsection{Proof of Theorem~\ref{t:main}}
To prove the convergence
of the rescaled clock process 
$\bar S_N(\cdot)=e^{-\gamma N} S_N(\cdot r(N))$ to the
stable subordinator $V_{\gamma /\beta^2}$, we 
check first the convergence of
finite-dimensional marginals.  
As can be guessed,  Proposition~\ref{p:comp} will serve to this purpose.  
Let $\ell$, $\{u_i\}$ and $\{t_i\}$ be as above. Then,
\begin{equation}
  \begin{split}
    \mathbb E\Big[&
      \exp\Big\{-\sum_{i=1}^\ell 
        u_i\big(\bar S_N(t_k)-\bar S_N(t_{k-1})\big)\Big\}
      \Big|Y_N\Big]
    \\&=
    \mathbb E\big[F_N(X_N^0;\{t_i\},\{u_i\})\big|Y_N\big]=
    \mathbb E\big[F_N(X_N^1;\{t_i\},\{u_i\})\big]+o(1),
  \end{split}
\end{equation}
as follows from Proposition~\ref{p:comp}.

The value of $ \mathbb E\big[F_N(X_N^1;\{t_i\},\{u_i\})\big]$ is not
difficult  to calculate. Define $j_N(i)=\lfloor t_i r(N)/\nu \rfloor$. Then
\begin{equation}
  \begin{split}
    \label{e:pra}
    \mathbb E&\big[F_N(X_N^1;\{t_i\},\{u_i\})\big]=
    \mathbb E\Big[\exp\Big(
        -\sum_{k=1}^\ell
        \frac {u_k} {e^{\gamma N}} 
        \sum_{i=t_{k-1}r(N)}^{t_k r(N)-1} e_i
        e^{\beta \sqrt N X^1_N(i)} \Big)
      \Big]
    \\&\ge
    \mathbb E\Big[
      \prod_{k=1}^\ell
      \prod_{j=j(k-1)+1}^{j(k)}
      \exp\Big(-
        \frac {u_k} {e^{\gamma N}} 
        \sum_{i=0}^{\nu-1} e_{j\nu +i}
        e^{\beta \sqrt N X^1_N(j\nu +i)} \Big)
      \Big]
  \end{split}
\end{equation}
Since the process  $X_N^1$ is a piece-wise independent process, the
product in \eqref{e:pra} is a product of independent random variables.
Then expectations of all of them can be then bounded 
using Proposition~\ref{p:oneblock}. We get, for $\delta >0$ fixed and 
$N$ large enough,
\begin{equation}
  \begin{split}
    \mathbb E\big[F_N(X_N^1&;\{t_i\},\{u_i\})\big]
    \ge
    \prod_{k=1}^\ell
    \prod_{j=j_N(k-1)+1}^{j_N(k)}
    \mathcal F_N(u_k ) 
    \\&\ge
    \prod_{k=1}^\ell
    \big(1-(1+\delta ) 
      \nu N^{-1/2}
      e^{-N\gamma^2/2\beta^2}Ku_k^{\gamma /\beta^2}\big)
    ^{j_N(k)-j_N(k-1)-1}
    \\&\ge
    \prod_{k=1}^\ell
    \exp\big\{-(1+2\delta )(t_k-t_{k-1})K u^{\gamma /\beta^2}\big\},
  \end{split}
\end{equation}
which is (up to $1+2\delta $ term) the Laplace transform of 
$V_{\gamma /\beta^2}(K\cdot)$. A corresponding upper bound can be
constructed analogously.

To check the tightness for $\bar S_N$ in $D([0,T],\mathbb R)$ equipped 
with the Skorokhod $M_1$-topology we use Theorem~\ref{t:tight}. Since 
the  processes $\bar S_N$ are increasing, it is easy to see that condition (i) 
is equivalent to the tightness of the distribution of $\bar S_N(T)$, 
which can be checked easily from the convergence of the Laplace transform 
of the marginal at time $T$ (the limiting 
Laplace transform tends to 1 as $u\to 0$).  

In order to check condition (ii), remark that for increasing functions the 
oscillation function $w'_{\bar S_N}(\delta )$ is always equal to zero. So 
checking (ii) boils down to controlling the boundary oscillations 
$v_{\bar S_N}(0,\delta )$ and $v_{\bar S_N}(T,\delta )$. 
For the first quantity (using again the monotonicity of $\bar S_N$) 
this amounts to check that 
$\mathbb P[\bar S_N(\delta )\ge \eta ]<\varepsilon $ if $\delta $ is
small enough and  $N$ large enough. Using  
the convergence of of marginal at time $\delta $, it is sufficient to
take $\delta$ such that 
$\mathbb P[V_{\gamma /\beta^2}(K\delta )\ge \eta ]\le \varepsilon /2$, and
take $n_0$ such that for all $n\ge n_0$ 
\begin{equation}
  \big |\mathbb P[\bar S_N(\delta )\ge \eta ]-
  \mathbb P[V_{\gamma /\beta^2}(K\delta )\ge \eta ]\big|
  \le \varepsilon /2.
\end{equation}
The reasoning for $v_{\bar S_N}(T,\delta )$ 
is analogous. 
\qed

\subsection{Coarse-grained clock process}
To prove our aging result, that is Theorem~\ref{t:aging}, we need to modify the result of
Theorem~\ref{t:main} slightly.
Let 
$\tilde S_N$  be the ``coarse-grained'' clock processes,
\begin{equation}
  \tilde S_N(t)=
  \frac 1 {e^{\gamma N}}
    S_N(\nu \lfloor t r(N)\nu^{-1}\rfloor ).
\end{equation}
For these processes we can strengthen the topology used in
Theorem~\ref{t:main}, that is we can  replace the $M_1$-
by the $J_1$-topology.

\begin{theorem}
  \label{t:J}
  If the hypothesis of Theorem~\ref{t:main} is satisfied, then
  \begin{equation}
    \tilde S_N(t)
    \xrightarrow{N\to\infty}
    V_{\gamma /\beta^2}(K t)
    \qquad \mathcal Y-\text{a.s.,}
  \end{equation}
  weakly in the $J_1$-topology on the space of càdlàg
  functions $D([0,T],\mathbb R)$.
\end{theorem}

Unfortunately, we cannot prove the theorem with estimates we have already
at disposition. We should return back and improve some of them. First we
show that traps with energies ``much smaller'' than 
$\gamma \sqrt N / \beta $ almost do not contribute to the clock process.
Let $B_m=\gamma \sqrt N /\beta - m /(\beta \sqrt N)$ and let 
\begin{equation}
  \bar S^m_N(t)=
  e^{-\gamma N}
  \sum_{i=0}^{\lfloor t r(N) \rfloor} 
  e_i \exp\big\{\beta \sqrt N X_N^0(i)\big\} 
  \bbone\{X_N^0(i)\le B_m\}.
\end{equation}

\begin{lemma}
  \label{l:small}
  For every $T$ and $\eta$, $\varepsilon   >0$ there exists $m$ large 
  enough such that
  \begin{equation}
    \mathbb P[\bar S^m_N(T)\ge \eta |\mathcal Y] \le \varepsilon ,
    \qquad \mathcal Y\text{-a.s.}
  \end{equation}
\end{lemma}
\begin{proof}
  To prove this lemma we should improve/modify slightly the calculations of
  Sections~\ref{s:oneblock} and~\ref{s:compar}. With the notation of
  Section~\ref{s:oneblock} define
  \begin{equation}
    \label{e:sc}
    \mathcal F^m_N =
    \mathbb E\Big[\exp \Big\{
        -e^{-\gamma N} \sum_{i=1}^\nu e_i e^{\beta \sqrt N U_i}
        \bbone\{U_i\le B_m\}
        \Big\}\Big].
  \end{equation}
  (comparing with \eqref{e:mathcalf} observe that we set $u=1$).
  We will show that
  \begin{equation}
    \label{e:sa}
    \lim_{N\to\infty}
    f(N) 
    e^{N \gamma^2/2\beta^2}\,
    [ 1-\mathcal F^m_N]=
    K_m,
  \end{equation}
  with $K_m\to 0$ as $m\to \infty$. The proof of this claim is completely
  analogous to the proof of Proposition~\ref{p:oneblock}. One should only
  modify the domains of integrations. More precisely, the definition of 
  $D_k$ which appears after \eqref{e:aa} should be replaced by
  $D_k^m=D_k\cap \{z:G_k(z)\le B_m\}$. Hence, $D'_k$
  becomes $D'^m_k=D'_k \cap\{b:G_k(b)\le -m/(\beta /\sqrt N)\}$,
  which then restricts the domain of integration in \eqref{e:vv} to 
  $(-\infty, -m/\beta ]$. Hence, the constant $K_m$ can be made arbitrarily
  small by choosing $m$ large. 

  Further, as in Section~\ref{s:compar}, define 
  \begin{equation}
    \begin{split}
      F_N^m(X)
      =
      \exp\Big(- 
        \sum_{i=0}^{T r(N)-1} 
        g\Big(e^{-\gamma N} 
          e^{\beta \sqrt N X(i)}\bbone\{X(i)\le B_m\}
          \Big)
        \Big).
    \end{split}
  \end{equation}
  Then, as in Proposition~\ref{p:comp}, we will show
  \begin{equation}
    \label{e:sd}
    \lim_{N\to \infty} 
    \mathbb E\big[F_N^m(X^0_N)\big|\mathcal Y\big]
    -\mathbb E\big[F_N^m(X^1_N)\big]=0,\qquad
    \text{$\mathcal Y$-a.s.}
  \end{equation}
  We use again \eqref{e:comparison} to show this claim.
  Although the indicator function is not differentiable, we will proceed as
  if it was, setting $(\bbone \{x\le B\})'=-\delta (x-M)$, where 
  $\delta $ denotes the Dirac delta function. As usual, this can be
  justified e.g.~by using smooth approximations of the indicator function. 
  The second derivative of $F_N^m(X)$ equals
  \begin{equation}
    \begin{split}
      \label{e:sb}
      &\frac{u^2 \beta^2 N} {e^{2\gamma N}} 
      e^{\beta \sqrt N(X(i)+X(j))}
      g'\big(u  
        e^{\beta \sqrt N X(i)-\gamma N}\big)  
      g'\big(u  
        e^{\beta \sqrt N X(j)-\gamma N}\big)     
      F^m_N(X)
      \\&\qquad\times 
      \Big(\bbone\{X(i)\le B_m\}-
        \frac{\delta_{B_m}(X(i))}{\beta \sqrt N}\Big)
      \Big(\bbone\{X(j)\le B_m\}-
        \frac{\delta_{B_m}(X(j))}{\beta \sqrt N}\Big)
      \\&\le 
      u^2 \beta^2 N 
      e^{\beta \sqrt N(X^h_N(i)+X^h_N(j))-2\gamma N}
      \exp\big(
        -2g\big(u
          e^{\beta \sqrt N X^h_N(i)-\gamma N}\big)
        -2g\big(u
          e^{\beta \sqrt N X^h_N(j)-\gamma N}\big)
        \big)
      \\&\qquad\times 
      \Big(\bbone\{X(i)\le B_m\}-
        \frac{\delta_{B_m}(X(i))}{\beta \sqrt N}\Big)
      \Big(\bbone\{X(j)\le B_m\}-
        \frac{\delta_{B_m}(X(j))}{\beta \sqrt N}\Big).
    \end{split}
  \end{equation}
  We should now bound the contributions of four terms. The one with the 
  product of two indicator functions is easy, because we can use directly 
  the result of Lemma~\ref{l:tech}. For remaining three terms, those with 
  the product of one indicator and one delta function, and this with two 
  delta function, the calculation should be repeated. However, in the end 
  we find that \eqref{e:sb} is bounded by $\bar \Xi (\Cov(X(i),X(j)))$ 
  as before. The presence of the delta functions makes actually the 
  calculations slightly less complicated. The proof then proceed as in 
  Section~\ref{s:compar}. 

  We can now finish the proof of Lemma~\ref{l:small}. By \eqref{e:sc} and
  \eqref{e:sd},
  \begin{equation}
    \begin{split}
      \mathbb E&\big[\exp(-\bar S^m_N(T))\big|\mathcal Y\big]=
      \mathbb E\big[F_N^m(X^0_N)\big|\mathcal Y\big]=
      \mathbb E\big[F_N^m(X^1_N)\big|\mathcal Y\big]+o(1)
      \\&=
      (1-K_m f(N)^{-1} e^{-N \gamma^2/2\beta^2})^{Tr(N)/\nu }+o(1)
      =
      e^{-K_m T} + o(1).
    \end{split}
  \end{equation}
  Since $K_m\to 0$ as $m\to \infty$,
  \begin{equation}
    \mathbb P[\bar S^m_N(T)\ge \eta |\mathcal Y]\le
    \frac{1-\mathbb E\big[\exp(-\bar S^m_N(T))\big|\mathcal Y\big]}
    {1-e^{-\eta }}
  \end{equation}
  can be made arbitrarily small by taking $m$ large enough.
\end{proof}

We study now how the blocks where the process visits sites with
energies larger than $B_m$ are distributed along the trajectory. To this
end we set for any Gaussian process $X$
\begin{equation}
  s_N^m(i;X)=\bbone\{\exists j: i\nu < j \le (i+1)\nu ,X(j)>B_m\}.
\end{equation}
and we define point process $H^m_N(X)$ on $[0,T]$ by
\begin{equation}
  H_N^m(X;\d x) = \sum_{i=0}^{T r(N)/\nu } s_N^m(i;X)\delta_{i\nu /r(N)}(\d x).
\end{equation}
\begin{lemma}
  \label{l:poinprocess}
  For every $m\in \mathbb R$ the point processes $H_N^m(X^0_N)$ converge to a
  homogeneous Poisson point process on $[0,T]$ with intensity 
  $\rho_m\in (0,\infty)$, $\mathcal Y$-a.s.
\end{lemma}
\begin{proof}
  To show this lemma we use Proposition 16.17 of Kallenberg \cite{Kal02}. 
  According to it, to prove the convergence of $H_N^m(X^0_N)$ to a Poisson 
  point process with intensity $\rho_m$ it is sufficient to check that for 
  any interval $I\subset[0,T]$
  \begin{equation}
    \label{e:kala}
    \lim_{N\to\infty} 
    \mathbb P[H_N^m(X^0_N;I)=0|\mathcal Y]=e^{-\rho_m |I|}
  \end{equation}
  and 
  \begin{equation}
    \label{e:kalb}
    \limsup_{N\to\infty}
    \mathbb E[H_N^m(X^0_N;I)|\mathcal Y]\le\rho_m|I|,
  \end{equation}
  where $|I|$ denotes the Lebesgue measure of $I$.
 
  The proof of the first claim is completely similar to the previous ones. We 
  start with a one-block estimate for \eqref{e:kala}:
  \begin{equation}
    \label{e:ba}
    \lim_{N\to\infty}
    N^{1/2}\nu^{-1} 
    e^{N \gamma^2/2\beta^2}\,
    \mathbb  E[s_N^m(0,U)]= \rho_m,
  \end{equation}
  Using the notation of Section~\ref{s:oneblock}, we get
  \begin{equation}
    \mathbb E[s_N^m(0,U)]=
    \int_{A_m}\frac {\d z}{(2\pi )^{\nu /2}}
    e^{-\frac 12 \sum_{i=1}^\nu z_i^2},
  \end{equation}
  where $A_m=\{z:\exists k \in \{1,\dots,\nu \} G_k(z)>B_m\}$. Dividing the
  domain of integration according to the maximal $G_k(z)$, this is equal
  \begin{equation}
    \sum_{k=1}^\nu 
    \int_{D_k}\frac {\d z}{(2\pi )^{\nu /2}}
    e^{-\frac 12 \sum_{i=1}^\nu z_i^2},
  \end{equation}
  where $D_k=\{z: G_k(z)>B_m, G_i(z)\le G_k(z)\forall i\neq k\}$. Using the
  substitution $z_i=b_i\pm\Gamma_i B_m$ on $D_k$ (where $+$ sign is used
    for $i\le k$ and $-$ sign for $i>k$) we get
  \begin{equation}
    e^{-N\gamma^2/2\beta^2}
    e^{m\gamma /\beta^2}\sum_{k=1}^\nu 
    \int_{D'_k}
    \frac {\d b}{(2\pi )^{\nu /2}}
    e^{-\frac 12 \sum_{i=1}^\nu b_i^2}e^{-B_mG_k(b)},
  \end{equation}
  where 
  $D'_k=\{b:G_k(b)>0, \sum_{j=i+1}^k b_j+ |k-i| \Gamma_\nu B_m\ge 0\forall i\neq k\}$. 
  The same reasoning as before then allows to show that the last 
  expression behaves like $\rho_m \nu  N^{-1/2}e^{-\gamma^2N/2\beta^2}$ as 
  $N\to \infty$.

  To compare the real process with the block-independent process, let
  \begin{equation}
    F_N(I;X)=\bbone\{\max\{X(i):i\nu /r(N)\in I\}\le B_m\}.
  \end{equation}
  The difference between $\mathbb E[F_N(I;X^0_N)|\mathcal Y]$ and 
  $\mathbb E[F_N(I;X^1_N)]$ is again given by the Gaussian comparison formula
  \eqref{e:comparison}. This time the second derivative equals
  \begin{equation}
    \delta (X(i)-B_m) \delta (X(j)-B_m) 
    \prod_{k\neq i,j} \bbone\{X(k)\le B_m\}
    \le
    \delta (X(i)-B_m) \delta (X(j)-B_m) .
  \end{equation}
  If covariance of $X(i)$ and $X(j)$ equals $c$, the expectation of the
  last expression is given by the value of the joint density of $X(i)$, 
  $X(j)$ at point $(B_m,B_m)$ which is
  \begin{equation}
    \label{e:bc}
    (2\pi (1-c^2))^{-1} e^{-B_m^2/(1+c)}\le
    C (1-c^2)^{-1}\exp\Big\{-\frac{\gamma^2 N}{\beta^2(1+c)}\Big\}.
  \end{equation}
  The exponential term is the same as in $\bar \Xi(c)$. The polynomial
  prefactor is however different, it diverges faster as $c\to 1$.
  We should thus return to \eqref{e:bigest} with $\tilde \Xi $ replaced by
  the right-hand side of \eqref{e:bc}. First
  \begin{equation}
    \int_0^1 (1-c^2)^{-1}=c^{-1}\arg\tanh(c)\approx -\frac 12 \log (1-c) 
  \end{equation}
  as $c\to 1$, which is not bounded for all $c$ as before. The estimates 
  \eqref{e:terma} and \eqref{e:termb} are influenced by this change. For 
  \eqref{e:terma} we can actually neglect this change, because the main 
  contribution to this term came from the neighborhood of $d=N/2$ (or 
    $c=0$) and was exponentially small  in the neighborhood of $d=1$ 
  (or $c\sim 1/N$). In the treatment of \eqref{e:termb}, the change has
  more effect, after some computations \eqref{e:termbint} becomes
  \begin{equation}
    Ct N^{3/2} \nu^{-1} \int_0^{\delta '}
    \log(c/x) e^{-c N x} \d x \le C t N^{1/2}\nu^{-1}\log N
    \xrightarrow{N\to\infty} 0.
  \end{equation}
  Finally, the change of polynomial prefactor of $\bar \Xi$ implies change
  in the control of \eqref{e:termc}. The equation~\eqref{e:termcsum}
  becomes
  \begin{equation}
    \eqref{e:termc} \le 
    C\sum_{d=0}^\nu 
    t N^{-3/2} d^2 
    [1-(1-2dN^{-1})^2p]^{-1}
    \exp(N\tilde\Upsilon (d/N)).
  \end{equation}
  and the linearization 
  of $\tilde\Upsilon $ gives new form of \eqref{e:termcint} 
  \begin{equation}
    C t N^{3/2}  \int_0^\varepsilon x e^{-c'N x}\d x 
    \le 
    C t N^{-1/2}
    \xrightarrow{N\to\infty} 0.
  \end{equation}
  Therefore, using \eqref{e:ba}
  \begin{equation}
    \begin{split}
      \mathbb P[H_N^m(X^0_N;I)=0|\mathcal Y]&=
      \mathbb E[F_N(I;X^0_N)|\mathcal Y]=
      \mathbb E[F_N(I;X^1_N)] +o(1)
      \\&=
      (1-\mathbb E[s^m_N(0,U)])^{|I|r(N)/\nu }
      \to e^{-\rho_m|I|}.
    \end{split}
  \end{equation}
  This completes the proof of \eqref{e:kala}.

  It is easy to check \eqref{e:kalb}. By definition,
  \begin{equation}
    \label{e:kalba}
    \mathbb E[H_N^m(X^0_N;I)|\mathcal Y]=
    \sum_{i:i\nu /R\in I} \mathbb E[s^m_N(i,X^0_N)|\mathcal Y].
  \end{equation}
  Since $\Lambda^0_{ij}\ge \Lambda^1_{ij}$ for $i$, $j$ in the same block, 
  $ \mathbb E[s^m_N(i,X^0_N)|\mathcal Y]\le \mathbb E[s^m_N(i,X^1_N)]$. Therefore,
  \begin{equation}
    \eqref{e:kalba}\le |I|r(N)/\nu \mathbb E[s^m_N(0,U)] = \rho_m|I|.
  \end{equation}
  This completes the proof of Lemma~\ref{l:poinprocess}.
\end{proof}

\begin{proof}[Proof of Theorem~\ref{t:J}]
  Checking the convergence of finite-dimensional marginals as well of
  condition (i) and the second part of (ii) of Theorem~\ref{t:tight} is
  analogous as for the original clock process $\bar S_N$. We should thus
  only prove the first part of condition~(ii). Namely that, for any 
  $\eta $ and $\varepsilon $ there exist $\delta $ such that
  \begin{equation}
    \label{e:Ja}
    \mathbb P[w_{\bar S_N} (\delta )\ge \eta ]\le \varepsilon, 
  \end{equation}
  for all $N$ large enough.
  
  Let 
  \begin{equation}
    w_f([\tau ,\tau +\delta ]) =
    \sup\{\min(|f(t_2)-f(t)|,|f(t)-f(t_1)|): 
      \tau \le t_1\le t \le t_2 \le \tau +\delta \}.
  \end{equation} 
  Fix $m$ such that 
  $\mathbb P[\bar S^m_N(T)\ge \eta /2]\le \varepsilon /2$, which is
  possible according to Lemma~\ref{l:small}. 
  If $H^m_N(X^0_n;[\tau ,\tau +\delta ])\le 1$ then 
  \begin{equation}
    w_{\bar S_N}([\tau ,\tau +\delta ])
    \le 
    \bar S^m_N(\tau +\delta )-\bar S^m_N(\tau )
    \le \bar S^m_N(T).
  \end{equation}
  Hence,
  \begin{equation}
    \mathbb P[ w_{\bar S_N}([\tau ,\tau +\delta ])\ge \eta |i
      \bar S^m_N(T)\le \eta /2  ]
    \le 
    \mathbb P[H^m_N(X^0_N;[\tau ,\tau +\delta ])\ge 2]
    \le C \rho_m \delta^2.
  \end{equation}

  We can now show \eqref{e:Ja}. Estimate
  \begin{equation}
    w_{\tilde S_N}(\delta )\le 
    \max \{ w_{\tilde S_N}([\tau ,\tau +2\delta ]):0\le \tau \le T, \tau =k
      \delta , k\in \mathbb N\}
  \end{equation}
  yields
  \begin{equation}
    \begin{split}
      \mathbb P[w_{\tilde S_N}(\delta )\ge \eta|\mathcal Y ]
      &\le 
      \sum_{k=0}^{T\delta^{-1}}
      \mathbb P[w_{\tilde S_N}([k\delta ,(k+2)\delta ]) \ge \varepsilon
        |\mathcal Y]
      \\&\le 
      \mathbb P[\bar S^m_N(T)\ge \eta /2]+
      \sum_{k=0}^{T\delta^{-1}}
      \mathbb P[H^m_N(X^0_N;[k\delta  ,(k +2)\delta ])\ge 2]
      \\&\le 
      \varepsilon /2 + C T\delta^{-1}\rho_m \delta^2 \le \varepsilon 
    \end{split}
  \end{equation}
  if $\delta $ is chosen small enough. This completes the proof.
\end{proof}

\begin{proof}[Proof of Theorem~\ref{t:aging}]
  Let $\mathcal R_N$ be the range of the coarse grained process 
  $\tilde S_N$. Obviously, for any $1>\varepsilon >0$,
  \begin{equation}
    A^\varepsilon_N(t,s)\supset \{\mathcal R_N\cap (t,s) =\emptyset\},
  \end{equation}
  because if the above intersection is empty, then $\sigma_N$ makes less
  than $\nu $ steps in time interval $[te^{\gamma N},se^{\gamma N}]$, and
  thus the overlap of $\sigma_N(te^{\gamma N})$ and 
  $\sigma_N(se^{\gamma N})$ is  $O(\nu /N)$.

  If $\mathcal R_N\cap (t,s) \neq\emptyset$, than there exist
  $u$ such that $\tilde S_N(u)\in (t,s)$. Moreover, it follows from
  Theorem~\ref{t:J} that for any $\delta $ there exist $\eta $ such than
  \begin{equation}
    \mathbb P[\tilde S_N(u+\eta )\in (s,t)]\ge 1-\delta .
  \end{equation}
  This however means that the process $\sigma_N$ make at least $\eta r(N)$
  steps between times $t$ and $s$ and thus the overlap between 
  $\sigma_N(te^{\gamma N})$ and $\sigma_N(se^{\gamma N})$ is with high
  probability close to 0. 

  Hence $\mathbb P[A_N^\varepsilon (t,s)|\mathcal Y]$ 
  is very well approximated by 
  $\mathbb P[\mathcal R_N\cap (t,s) =\emptyset|\mathcal Y]$. Since stable
  subordinator does not hit points, that is 
  $\mathbb P[\exists u: V_{\gamma /\beta^2}(u)=t]=0$, and $\tilde S_N$
  converge in $J_1$-topology,
  \begin{equation}
    \mathbb P[\mathcal R_N\cap (t,s) =\emptyset|\mathcal Y]
    \xrightarrow{N\to\infty}
    \mathbb P[\{V_{\gamma /\beta^2}(u):u\ge 0\} \cap (s,t)=\emptyset],
  \end{equation}
  which, as follows from the arc-sine law for stable subordinators, is
  given by the formula \eqref{e:aging}.
\end{proof}
  

\def\cprime{$'$}

\end{document}